\documentstyle{amsppt}
\magnification1200
\binoppenalty=10000\relpenalty=10000\relax
\TagsOnRight
\loadbold
\font\cyrms=wncyr8 \relax
\font\cyris=wncyi8 \relax
\font\cyrm=wncyr10 \relax
\font\cyrb=wncyb10 \relax
\font\cyri=wncyi10 \relax
\font\cyrc=wncysc10 \relax
\font\cyrbx=wncyb10 scaled1200 \relax
\redefine\cite#1{\rom{[\cyrm#1]}}
\let\hlop\!
\let\phi\varphi
\let\epsilon\varepsilon
\let\kappa\varkappa
\let\le\leqslant
\let\ge\geqslant
\let\wt\widetilde
\let\limsup\varlimsup
\let\<\langle
\let\>\rangle
\define\ZZ{\Bbb Z}
\define\NN{\Bbb N}
\define\QQ{\Bbb Q}
\define\KK{\Bbb K}
\define\CC{\Bbb C}
\define\GG{\bold G}
\define\sA{\Cal A}
\define\0{\bold0}
\define\bn{\boldkey n}
\define\balpha{\boldsymbol\alpha}
\define\bbeta{\boldsymbol\beta}
\define\den{\operatorname{den}}
\define\diag{\operatorname{diag}}
\define\Tr{\operatorname{Tr}}
\define\trans#1{{}^t\!#1}
\topmatter
\title
\cyrbx Sokrawenie faktorialov
\endtitle
\author
\cyrc V.\,V.~Zudilin
\endauthor
\address\cyrms
Moskovski\char"1A{} gosudarstvenny\char"1A{}
universitet im.~M.\,V.~Lomonosova
\endaddress
\email
{\tt wadim\@ipsun.ras.ru}
\endemail
\abstract\nofrills{}\cyrms

V rabote izuqaet\hbox{s}ya arifmetiqeskoe svo\char"1Astvo,
pozvolyayuwee usilivat\char"7E{} nekotorye teoretiko-qislovye ocenki.
Izvestnye ranee rezul\char"7Etaty nosi\-li, kak pravilo, kaqestvenny\char"1A{} kharakter.
Prilozhenie poluqennykh v rabo\-te koliqestvennykh rezul\char"7Etatov
k klassu obobwennykh gipergeometriqes\-kikh
$G$-funkci\char"1A{} rasxiryaet mnozhestvo irracional\char"7Enykh qisel,
yavlyayuwikh\-sya znaqeniyami e1tikh funkci\char"1A.

Bibliografiya: 19 nazvani\char"1A.
\endabstract
\endtopmatter
\footnote""{\cyrms Rabota vypolnena pri qastiqno\char"1A{} podderzhke Rossi\char"1Askogo
fonda fundamental\char"7Enykh issledovani\char"1A{} (grant~97-01-00181)
i ob\char"7Fedinennogo proekta fonda INTAS i
Rossi\char"1A\-skogo fonda fundamental\char"7Enykh issledovani\char"1A{}
(grant~{\rm IR-97-1904}).}
\leftheadtext{\cyrms V.\,V.~Zudilin}
\rightheadtext{\cyrms Sokrawenie faktorialov}
\document

\subhead\cyrb
1. Vvedenie
\endsubhead\cyrm
Vynesennoe v nazvanie stat\char"7Ei arifmetiqeskoe svo\char"1Astvo
ne obladaet mnogovekovo\char"1A{} istorie\char"1A, khotya nekotorye ego
proyavleniya otnosyat\-\hbox{s}ya k razryadu xkol\char"7Enykh zadaq.

\example{\cyrc\indent Primer 1}\cyrm
Operator~$D=d/dz$ differencirovaniya po peremenno\char"1A~$z$
perevodit kol\char"7Eco mnogoqlenov $\ZZ[z]$ s celymi koe1fficientami
v sebya; sle\-dovatel\char"7Eno, e1tim zhe svo\char"1Astvom obladayut i operatory~$D^n$,
$n=0,1,2,\dots$\,. Kak neslozhno pokazat\char"7E{} (sm., naprimer,
\cite{1, gl.~4, lemma~7}), pod de\char"1Astviem operatorov $D^n/n!$
kol\char"7Eco~$\ZZ[z]$ takzhe perekhodit v sebya:
$$
\frac1{n!}D^n=\frac1{n!}\frac{d^n}{dz^n}\:\ZZ[z]\to\ZZ[z],
\qquad n=0,1,2,\dots\,.
$$
\endexample\cyrm

\example{\cyrc\indent Primer 2}\cyrm
Posledovatel\char"7Enost\char"7E{} mnogoqlenov
$$
\<\lambda\>_0=1, \qquad
\<\lambda\>_n=\lambda(\lambda-1)\dotsb(\lambda-n+1),
\quad n=1,2,\dots,
\tag1
$$
soderzhit\hbox{s}ya v kol\char"7Ece $\ZZ[\lambda]$; znaqit, pri celykh~$\lambda$
mnogoqleny~\thetag{1} prinimayut celye znaqeniya. Na samom dele,
e1tim svo\char"1Astvom takzhe obladayut mnogoqleny
$$
\Delta_n(\lambda)=\frac{\<\lambda\>_n}{n!},
\qquad n=0,1,2,\dots
\tag2
$$
(sm., naprimer, \cite{2, otd.~8, gl.~2, zadaqa~84}), ne prinadlezhawie
pri $n\ge2$ kol\char"7Ecu $\ZZ[\lambda]$. Mnogoqleny~\thetag{2} nazyvayut\hbox{s}ya
{\cyri celoznaqnymi}.
\endexample\cyrm

Pust\char"7E{} pole $\KK$~-- nekotoroe algebraiqeskoe rasxirenie polya~$\QQ$
raciona\-l\char"7Enykh qisel, $\ZZ_{\KK}$~-- kol\char"7Eco celykh polya~$\KK$.
Vmesto kol\char"7Eca~$\ZZ[z]$ v primere~1 mozhno rassmatrivat\char"7E{}
kol\char"7Eco~$\ZZ_{\KK}[z]$:
$$
\frac1{n!}\frac{d^n}{dz^n}\:\ZZ_{\KK}[z]\to\ZZ_{\KK}[z],
\qquad n=0,1,2,\dots\,.
$$

\definition{\cyrc\indent Opredelenie 1}\cyrm
Budem govorit\char"7E, qto {\cyri operator $D\:\KK[z]\to\KK[z]$
udovlet\-voryaet usloviyu sokraweniya faktorialov s postoyanno\char"1A\/}
$\Psi\ge1$, esli suwestvuet posledovatel\char"7Enost\char"7E{} natural\char"7Enykh qisel
$\{\psi_k\}_{k\in\NN}$ takaya, qto
$$
\psi_k\frac{D^n}{n!}\:\ZZ_{\KK}[z]\to\ZZ_{\KK}[z],
\qquad n=0,1,\dots,k, \quad k\in\NN,
\qquad\text{i}\qquad
\limsup_{k\to\infty}\psi_k^{1/k}\le\Psi.
$$
\enddefinition\cyrm

Takim obrazom, spravedlivo sleduyuwee utverzhdenie.

\proclaim{\cyrc\indent Lemma 1}\cyri
Dlya lyubogo algebraiqeskogo rasxireniya~$\KK$ polya~$\QQ$
operator $d/dz\:\KK[z]\to\KK[z]$
udovletvoryaet usloviyu sokraweniya faktorialov s postoyanno\char"1A~$1$.
\endproclaim\cyrm

Vybiraya v primere~2 \ $\lambda=a/b\in\QQ$,
gde $a\in\ZZ$ i $b\in\NN$ vzaimno prosty,
zametim, qto obwi\char"1A{} znamenatel\char"7E{} qisel
$$
b^n\frac{\<\lambda\>_n}{n!}
=\frac{a(a-b)(a-2b)\dotsb(a-(n-1)b)}{n!},
\qquad n=1,2,\dots,k,
\tag3
$$
raven
$\prod_{p\mid b}p^{\tau_p(k)}$
(sm.~\cite{1, gl.~1, lemma~8} ili~\cite{3, gl.~\rom{I}, prilozhenie}), gde
$$
\tau_p(k)=\biggl[\frac kp\biggr]
+\biggl[\frac k{p^2}\biggr]+\biggl[\frac k{p^3}\biggr]+\dotsb
\le\frac k{p-1}
\tag4
$$
-- stepen\char"7E{} vkhozhdeniya prostogo qisla~$p$ v~$k!$\,.
Takim obrazom, obwi\char"1A{} znamenatel\char"7E{} qisel~\thetag{3} ne prevoskhodit
$e^{k\chi(b)}$, gde
$$
\chi(b)=\sum_{p\mid b}\frac{\log p}{p-1}, \qquad b\in\NN,
\tag5
$$
t.e\. imeet geometriqeski\char"1A{} poryadok rosta pri $k\to\infty$.

Razumeet\hbox{s}ya, v kaqestve oblasti opredeleniya mnogoqlenov~\thetag{1}
mozhno rassmatrivat\char"7E{} algebraiqeskoe rasxirenie~$\KK$ polya~$\QQ$
(ili dazhe kol\char"7Eco kvadratnykh matric s koe1fficientami iz~$\KK$,
sm\. dalee p.~4). Opredelim srazu {\cyri znamenatel\char"7E\/} $\den\lambda$
qisla $\lambda\in\KK$ kak naimen\char"7Exee natural\char"7Enoe~$b$ takoe,
qto $b\lambda\in\ZZ_{\KK}$.

\definition{\cyrc\indent Opredelenie 2}\cyrm
Budem govorit\char"7E, qto {\cyri e1lement~$\lambda$ polya~$\KK$
udovletvoryaet usloviyu sokraweniya faktorialov s postoyanno\char"1A\/}
$\Psi\ge1$, esli suwestvuet posledovatel\char"7Enost\char"7E{} natural\char"7Enykh qisel
$\{\psi_k\}_{k\in\NN}$ takaya, qto
$$
\psi_k\frac{\<\lambda\>_n}{n!}\in\ZZ_{\KK},
\qquad n=0,1,\dots,k, \quad k\in\NN,
\qquad\text{i}\qquad
\limsup_{k\to\infty}\psi_k^{1/k}\le\Psi,
$$
gde simvol $\<\,\cdot\,\>_n$ opredelyaet\hbox{s}ya
formulo\char"1A~\thetag{1}.
\enddefinition\cyrm

Podytozhim skazannoe k primeru~2 v vide sleduyuwego utverzhdeniya.

\proclaim{\cyrc\indent Lemma 2}\cyri
Racional\char"7Enoe qislo~$\lambda$
udovletvoryaet usloviyu sokraweniya faktorialov
s postoyanno\char"1A{} $be^{\chi(b)}$\rom, gde $b=\den\lambda$\rom,
a funkciya $\chi(\,\cdot\,)$ opredelyaet\hbox{s}ya formulo\char"1A~\thetag{5}.
\endproclaim\cyrm

\remark{\cyrc\indent Zameqanie}\cyrm
Racional\char"7Enost\char"7E{} qisla~$\lambda$ v lemme~2
suwestvenna. Rassmat\-rivaya pole $\KK=\QQ(\lambda)$
koneqno\char"1A{} stepeni $\kappa=[\KK:\QQ]\ge2$
i polagaya v osnovno\char"1A{} lemme raboty~\cite{4} \
$a_1=-\lambda-1$, $b_1=0$,
$m=m_1=1$, $m_2=0$, $\tau=1-1/\kappa$, $\epsilon=1/6$,
poluqaem sleduyuwee utverzhdenie:
{\cyri dlya lyubo\char"1A{} posledovatel\char"7Enosti natural\char"7Enykh qisel
$\{\psi_k\}_{k\in\NN}$ tako\char"1A\rom, qto
$$
\psi_k\frac{\<\lambda\>_n}{n!}\in\ZZ_{\KK},
\qquad n=0,1,\dots,k, \quad k\in\NN,
$$
vypolnyaet\hbox{s}ya neravenstvo
$$
\psi_k\ge Ck^{(\tau-\epsilon)\kappa k}\ge Ck^{2k/3},
\qquad k\in\NN,
\tag6
$$
gde polozhitel\char"7Enaya postoyannaya~$C$
zavisit tol\char"7Eko ot qisla\/}~$\lambda$.
Ocenka~\thetag{6} oznaqaet, qto v sluqae
irracional\char"7Enogo~$\lambda$ rost obwikh znamenatele\char"1A{}
po\-sledovatel\char"7Enosti~\thetag{2} pri $k\to\infty$
po kra\char"1Ane\char"1A{} mere
faktorial\char"7Eny\char"1A{} i govo\-rit\char"7E{} o sokrawenii
faktorialov ne imeet smysla.
\endremark\cyrm

\subhead\cyrb
2. Istoriya voprosa
\endsubhead\cyrm
Primery, privedennye v~p.~1, mozhno sqitat\char"7E{} klassiqeskimi.
Na samom dele, ponyatie {\cyri sokrawenie faktorialov\/}
poyavilos\char"7E{} v stat\char"7Ee A.\,I.~Galoqkina~\cite{5} i
bylo svyazano s sokraweniem koe1fficientov
line\char"1Anykh priblizhayuwikh form (priblizheni\char"1A{} Pade)
dlya tak nazyvaemykh $G$-funkci\char"1A{} v metode Zigelya--Xidlovskogo.
Namek na e1to obsto\hbox{ya}tel\char"7E\-stvo
soderzhit\hbox{s}ya v rabote K.~Zigelya~\cite{6},
s kotoro\char"1A{} i beret naqalo ukazanny\char"1A{} metod.
Vo vsekh rabotakh, svyazannykh s primeneniem
metoda Zigelya--Xid\-lovskogo k klassu
$G$-funkci\char"1A{} i opublikovannykh do poyavleniya \cite{5},
byli poluqeny ocenki line\char"1Anykh form i
mnogoqlenov ot znaqeni\char"1A{} e1tikh funkci\char"1A{} v racional\char"7Enykh
toqkakh, veliqina kotorykh zavisela ot vysoty
rassmat\-rivaemykh form (sm., naprimer,~\cite{7}). I imenno
uslovie sokrawe\-niya faktorialov pozvolilo poluqit\char"7E{} dlya
odnogo podklassa $G$-funkci\char"1A{} ocenki module\char"1A{} mnogoqlenov ot
ikh znaqeni\char"1A{} v toqkakh, veliqina kotorykh uzhe ne zavisela ot
vysoty mnogoqlenov. V~\cite{5} byla sformulirovana teorema ob
e1ffektivno\char"1A{} ocenke line\char"1Ano\char"1A{} formy ot $G$-funkci\char"1A{} iz odnogo
klassa, dokazatel\char"7Estvo kotoro\char"1A{} opublikovano v
rabote~\cite{8}. V 1985\,g\. D.\,V.~Qudnovski\char"1A{} i
G.\,V.~Qudnovski\char"1A~\cite{9} dokazali, qto uslovie
sokrawe\-niya faktorialov, sformulirovannoe v~\cite{5},
vypolnyaet\hbox{s}ya dlya odnorodnykh sistem line\char"1Anykh
differencial\char"7Enykh uravneni\char"1A, kotorym udovletvoryayut
$G$-funkcii. V~\cite{3, gl.~\rom{VI}, \S\,4} neslozhnoe obobwenie
konstrukcii Qudnovskikh pozvolilo rasprostranit\char"7E{} e1to uslovie
na neodnorodnye sistemy line\char"1Anykh differencial\char"7Enykh uravneni\char"1A.

Nesmotrya na polozhitel\char"7Enoe rexenie problemy sokraweniya faktorialov
dlya $G$-funkci\char"1A, postoyannaya, s kotoro\char"1A{} proiskhodit e1to sokrawenie,
``daleka ot soverxenstva''. E1to obstoyatel\char"7Estvo svyazano s tem,
qto konstrukciya Qudnovskikh ispol\char"7Ezovala postroeniya
s pomow\char"7Eyu principa Dirikhle, imeyuwego slixkom bol\char"7Exo\char"1A{} zapas proqnosti.
V nastoyawe\char"1A{} rabote privodit\hbox{s}ya novoe rexenie ukazanno\char"1A{} problemy
dlya obobwennykh gipergeometriqeskikh $G$-funkci\char"1A,
sformulirovannoe v stat\char"7Ee~\cite{10, \S\,12} v kaqestve gipotezy.
Naxe dokazatel\char"7Estvo opiraet\hbox{s}ya na toqnye postroeniya i poe1tomu daet
bolee toqnuyu ocenku postoyanno\char"1A, s kotoro\char"1A{} proiskhodit
sokrawenie faktorialov. Teoretiko-qislovoe primenenie
poluqennogo rezul\char"7Etata rasxiryaet mnozhestvo irracional\char"7Enykh qisel,
yavlyayuwikhsya
\linebreak
znaqeniyami gipergeometriqeskikh $G$-funkci\char"1A.

Opredelim klass $G$-funkci\char"1A, imeyuwikh racional\char"7Enye koe1fficienty
v razlozhenii Te\char"1Alora v okrestnosti nulya. Skazhem, qto
{\cyri sovokupnost\char"7E{} funkci\char"1A{}
$$
f_j(z)=\sum_{n=0}^\infty f_{jn}z^n, \quad j=1,\dots,m,
\qquad f_{jn}\in\QQ, \quad j=1,\dots,m, \quad n=0,1,\dots,
\tag7
$$
prinadlezhit klassu\/} $\GG(C,\Phi)$, esli funkcii~\thetag{7}
analitiqny v kruge $|z|<C$ i suwestvuet posledovatel\char"7Enost\char"7E{}
natural\char"7Enykh qisel $\{\phi_k\}_{k\in\NN}$ takaya, qto
$$
\phi_kf_{jn}\in\ZZ,
\qquad j=1,\dots,m, \quad n=0,1,\dots,k, \quad k\in\NN,
\quad\;\text{i}\;\quad
\limsup_{k\to\infty}\phi_k^{1/k}\le\Phi.
$$

Sformuliruem teper\char"7E{} uslovie sokraweniya faktorialov dlya sistem
line\char"1Anykh differencial\char"7Enykh uravneni\char"1A{}
$$
\gathered
\frac d{dz}y_l=Q_{l0}+\sum_{j=1}^mQ_{lj}y_j,
\qquad l=1,\dots,m,
\\
Q_{lj}=Q_{lj}(z)\in\QQ(z),
\qquad l=1,\dots,m, \quad j=0,\dots,m,
\endgathered
\tag8
$$
kotorym udovletvoryayut $G$-funkcii~\thetag{7}.

Pust\char"7E{} mnogoqlen $T(z)\in\QQ[z]$~-- znamenatel\char"7E{} racional\char"7Enykh
funkci\char"1A{} so starxim koe1fficientom ravnym~$1$:
$$
T(z)Q_{lj}(z)\in\QQ[z],
\qquad l=1,\dots,m, \quad j=0,\dots,m.
\tag9
$$

Iz~\thetag{8} sleduet, qto dlya proizvodnykh poryadka~$n$,
$n=1,2,\dots$, imeyut mesto sootnoxeniya
$$
\gathered
\frac{d^n}{dz^n}y_l=Q_{l0}^{[n]}+\sum_{j=1}^mQ_{lj}^{[n]}y_j,
\qquad l=1,\dots,m,
\\
Q_{lj}^{[n]}=Q_{lj}^{[n]}(z)\in\QQ(z),
\qquad l=1,\dots,m, \quad j=0,\dots,m.
\endgathered
\tag10
$$
Neslozhnye vykladki pokazyvayut spravedlivost\char"7E{}
sleduyuwikh rekurrentnykh sootnoxeni\char"1A:
$$
\gathered
Q_{lj}^{[n]}(z)
=\frac d{dz}Q_{lj}^{[n-1]}(z)
+\sum_{r=1}^mQ_{lr}^{[n-1]}(z)Q_{rj}(z),
\\
l=1,\dots,m, \quad j=0,1,\dots,m,
\quad n=1,2,\dots;
\endgathered
\tag11
$$
poe1tomu
$$
\gathered
\aligned
T^n(z)Q_{lj}^{[n]}(z)
&=T(z)\frac d{dz}\bigl(T^{n-1}(z)Q_{lj}^{[n-1]}(z)\bigr)
-(n-1)T'(z)\cdot T^{n-1}(z)Q_{lj}^{[n-1]}(z)
\\ &\qquad
+\sum_{r=1}^mT^{n-1}(z)Q_{lr}^{[n-1]}(z)\cdot T(z)Q_{rj}(z),
\endaligned
\\
l=1,\dots,m, \quad j=0,1,\dots,m,
\quad n=1,2,\dots\,.
\endgathered
$$
Ot\hbox{s}yuda, v qastnosti, sleduet, qto
$$
T^n(z)Q_{lj}^{[n]}(z)\in\QQ[z],
\qquad l=1,\dots,m, \quad j=0,\dots,m,
\quad n=1,2,\dots\,.
$$

\definition{\cyrc\indent Opredelenie 3}\cyrm
Budem govorit\char"7E, qto {\cyri sistema
line\char"1Anykh differencia\-l\char"7Enykh
uravneni\char"1A~\thetag{8} udovletvoryaet
usloviyu sokraweniya faktorialov
s postoyanno\char"1A\/} $\Psi\ge1$,
esli suwestvuyut natural\char"7Enye qisla
$\{\psi_k\}_{k\in\NN}$ takie, qto
$$
\gather
\psi_k\frac{T^n(z)Q_{lj}^{[n]}(z)}{n!}\in\ZZ[z],
\;\quad l=1,\dots,m, \;\; j=0,\dots,m,
\;\; n=0,1,\dots,k, \;\; k\in\NN,
\\
\limsup_{k\to\infty}\psi_k^{1/k}\le\Psi.
\endgather
$$
\enddefinition\cyrm

Uslovie sokraweniya faktorialov dlya sistemy~\thetag{8}
vypolnyaet\hbox{s}ya, voobwe govorya,
tol\char"7Eko esli e\char"1A{} udovletvoryaet
nabor line\char"1Ano nezavisimykh nad $\CC(z)$ \ $G$-funkci\char"1A.
Podobnye sistemy otnosyat\hbox{s}ya k klassu sistem diffe\-rencial\char"7Enykh
uravneni\char"1A{} fuksovskogo tipa.
S toqnost\char"7Eyu do meromorfnogo
preobrazovaniya prostranstva rexeni\char"1A{}
matrica koe1fficientov
$Q(z)=\allowmathbreak\bigl(Q_{lj}(z)\bigr)_{l,j}$
sistemy~\thetag{8} fuksovskogo tipa imeet vid
$$
Q(z)=\frac1{z-\gamma_1}A_1+\dots+\frac1{z-\gamma_s}A_s,
\tag12
$$
gde
$\gamma_1,\dots,\gamma_s$~-- regulyarnye osobennosti
sistemy~\thetag{8}, $A_1,\dots,A_s$~-- qislovye matricy
(sm\. \cite{11, zameqanie k~\S\,2.4}).
V sluqae~\thetag{12} znamenatel\char"7E{}
sootvet\-\hbox{s}tvuyuwe\char"1A{} sistemy~\thetag{8}
raven $T(z)=(z-\gamma_1)\dotsb(z-\gamma_s)$.

\subhead\cyrb
3. Sokrawenie faktorialov dlya differencial\char"7Enykh operatorov
\endsubhead\cyrm
V e1tom punkte my issleduem odno obobwenie
operatora differencirovaniya $d/dz$, dlya kotorogo
takzhe vypolneno uslovie sokraweniya faktorialov.

Pust\char"7E{} $\KK$~-- algebraiqeskoe rasxirenie polya~$\QQ$.
Rassmotrim differencial\char"7Eny\char"1A{} operator
$$
D=\frac d{dz}+\frac\lambda z,
\qquad \lambda\in\QQ.
\tag13
$$
Operator $T(z)D$, gde $T(z)=z$, perevodit kol\char"7Eco $\KK[z]$ v sebya,
v to vremya kak sam operator~$D$ otobrazhaet $\KK[z]$ v~$\KK(z)$.
Poe1tomu my neskol\char"7Eko rasxirim oblast\char"7E{} de\char"1Astviya opredeleniya~1.

\definition{\cyrc\indent Opredelenie 1$'$}\cyrm
{\cyri Znamenatelem\/} operatora $D\:\KK(z)\to\KK(z)$
my nazovem netrivial\char"7Eny\char"1A{} mnogoqlen $T(z)\in\KK[z]$
naimen\char"7Exe\char"1A{} stepeni so starxim koe1fficientom ravnym~$1$,
dlya kotorogo $T(z)D\:\KK[z]\to\KK[z]$.
Budem go\-vorit\char"7E, qto {\cyri operator~$D$
udovletvoryaet usloviyu sokraweniya faktorialov s postoyanno\char"1A\/}
$\Psi\ge1$, esli suwestvuet posledovatel\char"7Enost\char"7E{} natural\char"7Enykh qisel
$\{\psi_k\}_{k\in\NN}$ takaya, qto
$$
\psi_k\frac{T^n(z)D^n}{n!}\:\ZZ_{\KK}[z]\to\ZZ_{\KK}[z],
\;\quad n=0,1,\dots,k, \quad k\in\NN,
\;\quad\text{i}\;\quad
\limsup_{k\to\infty}\psi_k^{1/k}\le\Psi.
$$
\enddefinition\cyrm

\proclaim{\cyrc\indent Lemma 3}\cyri
Dlya differencial\char"7Enogo operatora~\thetag{13} spravedlivy
tozhdest\-va
$$
D^n=\sum_{l=0}^n\binom nl\frac{\<\lambda\>_l}{z^l}
\frac{d^{n-l}}{dz^{n-l}},
\qquad n\in\NN,
\tag14
$$
gde simvol $\<\,\cdot\,\>_l$ opredelyaet\hbox{s}ya formulo\char"1A~\thetag{1}.
\endproclaim\cyrm

\demo{\cyrc\indent Dokazatel\char"7Estvo}\cyrm
Pri $n=1$ formula~\thetag{14} sovpadaet s opredeleniem \thetag{13}
operatora~$D$. Polagaya, qto formula~\thetag{14}
spravedliva dlya nekotorogo natural\char"7Enogo~$n$, imeem
$$
\align
D^{n+1}
&=\biggl(\frac d{dz}+\frac\lambda z\biggr)D^n
=\frac d{dz}
\sum_{l=0}^n\binom nl\frac{\<\lambda\>_l}{z^l}
\frac{d^{n-l}}{dz^{n-l}}
+\frac\lambda z
\sum_{l=0}^n\binom nl\frac{\<\lambda\>_l}{z^l}
\frac{d^{n-l}}{dz^{n-l}}
\\
&=\sum_{l=0}^n\binom nl\frac{\<\lambda\>_l}{z^l}
\frac{d^{n-l+1}}{dz^{n-l+1}}
-\sum_{l=0}^n\binom nl\frac{l\cdot\<\lambda\>_l}{z^{l+1}}
\frac{d^{n-l}}{dz^{n-l}}
+\sum_{l=0}^n\binom nl\frac{\lambda\cdot\<\lambda\>_l}{z^{l+1}}
\frac{d^{n-l}}{dz^{n-l}}
\\
&=\sum_{l=0}^n\binom nl\frac{\<\lambda\>_l}{z^l}
\frac{d^{n+1-l}}{dz^{n+1-l}}
+\sum_{l=0}^n\binom nl\frac{\<\lambda\>_{l+1}}{z^{l+1}}
\frac{d^{n-l}}{dz^{n-l}}
\\
&=\sum_{l=0}^n\binom nl\frac{\<\lambda\>_l}{z^l}
\frac{d^{n+1-l}}{dz^{n+1-l}}
+\sum_{l=1}^{n+1}\binom n{l-1}\frac{\<\lambda\>_l}{z^l}
\frac{d^{n+1-l}}{dz^{n+1-l}}
\\
&=\sum_{l=0}^{n+1}\binom{n+1}l\frac{\<\lambda\>_l}{z^l}
\frac{d^{n+1-l}}{dz^{n+1-l}}.
\endalign
$$
Tem samym, formula~\thetag{14} spravedliva i dlya~$n+1$.
Soglasno principu matematiqesko\char"1A{} indukcii
ona verna dlya vsekh natural\char"7Enykh~$n$.
Lemma dokazana.
\enddemo\cyrm

Pol\char"7Ezuyas\char"7E{} opredeleniem binomial\char"7Enykh koe1fficientov,
iz lemmy~3 poluqaem tozhdestvo
$$
z^n\frac{D^n}{n!}
=\sum_{l=0}^nz^{n-l}\cdot\frac{\<\lambda\>_l}{l!}
\cdot\frac1{(n-l)!}\frac{d^{n-l}}{dz^{n-l}},
\qquad n\in\NN.
$$
Primenyaya teper\char"7E{} usloviya sokraweniya faktorialov dlya
operatora~$d/dz$ (lemma~1) i racional\char"7Enogo
qisla~$\lambda$ (lemma~2), poluqaem sleduyuwee utver\-zhdenie.

\proclaim{\cyrc\indent Teorema 1}\cyri
Differencial\char"7Eny\char"1A{} operator~\thetag{13} udovletvoryaet
usloviyu so\-kraweniya faktorialov s postoyanno\char"1A{}
$be^{\chi(b)}$\rom, gde $b=\den\lambda$\rom, a funkciya $\chi(\,\cdot\,)$
opredelyaet\hbox{s}ya formulo\char"1A~\thetag{5}.
\endproclaim\cyrm

Lemma~3 legko obobwaet\hbox{s}ya na sluqa\char"1A{}
differencial\char"7Enogo operatora
$$
D=\frac d{dz}
+\frac{\lambda_1}{z-\gamma_1}+\dots+\frac{\lambda_s}{z-\gamma_s},
\qquad
\lambda_1,\dots,\lambda_s,\gamma_1,\dots,\gamma_s\in\QQ,
\tag15
$$
so znamenatelem $T(z)=(z-\gamma_1)\dotsb(z-\gamma_s)$.
Ukazhem sootvet\hbox{s}tvuyuwie tozhdestva bez dokazatel\char"7Estva.

\proclaim{\cyrc\indent Lemma 4}\cyri
Dlya differencial\char"7Enogo operatora~\thetag{15}
spravedlivy tozhdest\-va
$$
D^n=\sum\Sb n_0,n_1,\dots,n_s\ge0\\n_0+n_1+\dots+n_s=n\endSb
\frac{n!}{n_0!\,n_1!\dotsb n_s!}
\frac{\<\lambda_1\>_{n_1}}{(z-\gamma_1)^{n_1}}
\dotsb
\frac{\<\lambda_s\>_{n_s}}{(z-\gamma_s)^{n_s}}
\frac{d^{n_0}}{dz^{n_0}},
\qquad n\in\NN.
\tag16
$$
\endproclaim\cyrm

Soglasno lemme~4 vypolneno
$$
\align
\frac{T^n(z)D^n}{n!}
=\sum\Sb n_0,n_1,\dots,n_s\ge0\\n_0+n_1+\dots+n_s=n\endSb
&
(z-\gamma_1)^{n-n_1}\dotsb(z-\gamma_s)^{n-n_s}
\\ &\quad\times
\frac{\<\lambda_1\>_{n_1}}{n_1!}\dotsb\frac{\<\lambda_s\>_{n_s}}{n_s!}
\cdot\frac1{n_0!}\frac{d^{n_0}}{dz^{n_0}},
\qquad n\in\NN,
\endalign
$$
a znaqit, differencial\char"7Eny\char"1A{} operator~\thetag{15}
takzhe udovletvoryaet usloviyu sokraweniya faktorialov.

\proclaim{\cyrc\indent Teorema 2}\cyri
Differencial\char"7Eny\char"1A{} operator~\thetag{15} udovletvoryaet
usloviyu so\-kraweniya faktorialov s postoyanno\char"1A{} $qbe^{\chi(b)}$\rom, gde
$q$~-- proizvedenie znamenatele\char"1A{} qisel $\gamma_1,\dots,\gamma_s$\rom,
$b=\den(\lambda_1,\dots,\lambda_s)$~--
naimen\char"7Exi\char"1A{} obwi\char"1A{} znamenatel\char"7E{}
qisel $\lambda_1,\dots,\lambda_s$\rom,
a funkciya $\chi(\,\cdot\,)$ opredelyaet\hbox{s}ya formulo\char"1A~\thetag{5}.
\endproclaim\cyrm

\subhead\cyrb
4. Sokrawenie faktorialov dlya kvadratnykh matric
\endsubhead\cyrm
Qislovuyu kvad\-ratnuyu matricu~$A$ razmera~$m$ mozhno podstavit\char"7E{}
v lyubo\char"1A{} mnogoqlen; v qastnosti,
$$
\<A\>_0=E, \qquad
\<A\>_n=A(A-E)\dotsb(A-(n-1)E),
\quad n=1,2,\dots,
$$
gde $E$~-- ediniqnaya matrica razmera~$m$.
Esli vse e1lementy matricy~$A$ prinadlezhat
algebraiqeskomu rasxireniyu~$\KK$ polya~$\QQ$, to
vpolne estest\-vennym yavlyaet\hbox{s}ya sleduyuwee opredelenie.

\definition{\cyrc\indent Opredelenie 2$'$}\cyrm
Skazhem, qto {\cyri matrica\/}~$A$
s e1lementami iz~$\KK$
{\cyri udovlet\-voryaet usloviyu sokraweniya faktorialov
s postoyanno\char"1A\/} $\Psi\ge1$,
esli suwestvuyut natural\char"7Enye qisla
$\{\psi_k\}_{k\in\NN}$ takie, qto e1lementy matric
$$
\psi_k\frac{\<A\>_n}{n!},
\qquad n=0,1,\dots,k, \quad k\in\NN,
$$
prinadlezhat~$\ZZ_{\KK}$ i
$$
\limsup_{k\to\infty}\psi_k^{1/k}\le\Psi.
$$
\enddefinition\cyrm

\proclaim{\cyrc\indent Lemma 5}\cyri
Esli matrica~$A$ udovletvoryaet usloviyu sokraweniya
faktorialov s postoyanno\char"1A~$\Psi$\rom, to tem zhe svo\char"1Astvom
obladaet matrica~$TAT^{-1}$\rom,
gde $T$~-- proizvol\char"7Enaya nevyrozhdennaya matrica\rom,
e1lementy kotoro\char"1A~-- algebraiqeskie qisla.
\endproclaim\cyrm

\demo\nofrills{\cyrc\indent Dokazatel\char"7Estvo}\cyrm\
osnovano na primenenii e1lementarnogo tozhdestva
$$
(TAT^{-1})^n=TA^nT^{-1}, \qquad n=0,1,2,\dots,
$$
otkuda, v qastnosti, sleduet, qto
$$
\frac{\<TAT^{-1}\>_n}{n!}=T\frac{\<A\>_n}{n!}T^{-1},
\qquad n=0,1,2,\dots\,.
$$
Esli $t_1,t_2$~-- naimen\char"7Exie obwie znamenateli e1lementov
matric~$T,T^{-1}$ so\-otvet\hbox{s}tvenno,
a $\{\psi_k\}_{k\in\NN}$~-- posledovatel\char"7Enost\char"7E{}
iz opredeleniya~2$'$ dlya mat\-ricy~$A$,
to v kaqestve sootvet\hbox{s}tvuyuwe\char"1A{}
posledovatel\char"7Enosti dlya mat\-ricy~$TAT^{-1}$
mozhno vzyat\char"7E{} $\{t_1t_2\psi_k\}_{k\in\NN}$.
Ostalos\char"7E{} zametit\char"7E, qto
$$
\limsup_{k\to\infty}(t_1t_2\psi_k)^{1/k}
=\limsup_{k\to\infty}\psi_k^{1/k}.
$$
Lemma dokazana.
\enddemo\cyrm

Iz lemmy~5 sleduet, qto sokrawenie faktorialov
i vyqislenie postoyanno\char"1A, s kotoro\char"1A{} ono proiskhodit,
dostatoqno dokazat\char"7E{} dlya matricy, imeyuwe\char"1A{} zhordanovu normal\char"7Enuyu formu.

\proclaim{\cyrc\indent Lemma 6}\cyri
Pust\char"7E{} matrica~$A$ sostoit iz zhordanovykh kletok $A_1,\dots,A_s$\rom,
raspolozhennykh vdol\char"7E{} glavno\char"1A{}
diagonali i otveqayuwikh \rom(ne obyazatel\char"7Eno
raz\-liqnym\rom) sobstvennym znaqeniyam $\lambda_1,\dots,\lambda_s$
sootvet\hbox{s}tvenno. Togda mat\-rica
$\Delta_n(A)=\<A\>_n/n!$ sostoit tol\char"7Eko iz kletok
$\Delta_n(A_1),\dots,\Delta_n(A_s)$\rom,
raspo\-lozhennykh vdol\char"7E{} glavno\char"1A{}
diagonali\rom, priqem
$$
\gathered
\Delta_n(A_l)=\pmatrix
\Delta_n(\lambda_j) & \frac1{1!}\Delta_n^{(1)}(\lambda_j)
& \frac1{2!}\Delta_n^{(2)}(\lambda_j) & \frac1{3!}\Delta_n^{(3)}(\lambda_j)
& \hdots \\
0 & \Delta_n(\lambda_j) & \frac1{1!}\Delta_n^{(1)}(\lambda_j)
& \frac1{2!}\Delta_n^{(2)}(\lambda_j) & \hdots \\
0 & 0 & \Delta_n(\lambda_j) & \frac1{1!}\Delta_n^{(1)}(\lambda_j)
& \hdots \\
\hdotsfor5 \\
0 & 0 & \hdots & 0 & \Delta_n(\lambda_j)
\endpmatrix,
\\
l=1,\dots,s.
\endgathered
$$
\endproclaim\cyrm

\demo{\cyrc\indent Dokazatel\char"7Estvo}\cyrm
E1to utverzhdenie khoroxo izvestno. Bolee togo,
ono ostaet\hbox{s}ya vernym, esli zamenit\char"7E{} mnogoqlen
$\Delta_n(\,\cdot\,)$ na lyubuyu druguyu funk\-ciyu
analitiqeskuyu v kruge $|z|<C$, soderzhawem sobstvennye znaqeniya
\linebreak
$\lambda_1,\dots,\lambda_s$
(sm., naprimer, \cite{12, gl.~5, \S\,1, primer~2}).
\enddemo\cyrm

\remark{\cyrc\indent Zameqanie 1}\cyrm
Kak sleduet iz lemmy~6, na glavno\char"1A{} diagonali
matricy $\Delta_n(A)$ stoyat qisla
$\Delta_n(\lambda_1),\dots,\Delta_n(\lambda_s)$;
v sootvet\hbox{s}tvii s zameqaniem k lem\-me~2 o sokrawenii
faktorialov dlya matricy~$A$ mozhno govorit\char"7E{} lix\char"7E{}
v sluqae racional\char"7Enykh $\lambda_1,\dots,\lambda_s$.
Poe1tomu v dal\char"7Ene\char"1Axem my ograniqimsya rassmotreniem
{\cyri racional\char"7Enykh kvadratnykh matric}, sobstvennye
znaqeniya kotorykh racional\char"7Eny. {\cyri Znamenatelem\/}
$\den A$ racional\char"7Eno\char"1A{} matricy~$A$ nazovem naimen\char"7Exi\char"1A{}
obwi\char"1A{} znamenatel\char"7E{} ee sobstvennykh znaqeni\char"1A.
Pri e1tom matrica~$bA$ ne obyazatel\char"7Eno imeet
celoqislennye e1lementy; primerom sluzhit matrica
$$
\pmatrix 1/2 & 1/2 \\ 1/2 & 1/2 \endpmatrix,
$$
znamenatel\char"7E{} kotoro\char"1A{} raven~$1$ (ee sobstvennye znaqeniya~$0,1$).
\endremark\cyrm

\remark{\cyrc\indent Zameqanie 2}\cyrm
Dlya kazhdogo sobstvennogo znaqeniya~$\lambda_l$
racional\char"7Eno\char"1A{} matricy~$A$ opredelim
natural\char"7Enoe qislo~$r_l$ kak maksimal\char"7Eny\char"1A{} razmer zhordanovo\char"1A{} kletki,
otveqayuwe\char"1A{} sobstvennomu znaqeniyu~$\lambda_l$,
v normal\char"7Eno\char"1A{} forme matricy~$A$. Otmetim, qto
$$
(\lambda-\lambda_1)^{r_1}
(\lambda-\lambda_2)^{r_2}\dotsb
(\lambda-\lambda_s)^{r_s},
\qquad \lambda_j\le\lambda_l, \quad j\ne l,
\tag17
$$
nazyvaet\hbox{s}ya {\cyri minimal\char"7Enym mnogoqlenom\/} matricy~$A$;
on delit lyubo\char"1A{} mnogoqlen $P(\lambda)$ tako\char"1A, qto $P(A)=0$.
V qastnosti, soglasno teoreme Gamil\char"7E\-tona--Ke1li
minimal\char"7Eny\char"1A{} mnogoqlen
delit kharakteristiqeski\char"1A{} mnogoqlen
$\det(A-\lambda E)$ matricy~$A$.
Krome togo, esli $\trans A$~-- matrica,
transponirovannaya k matrice~$A$,
to minimal\char"7Enye (kharakteristiqeskie)
mnogoqleny mat\-ric~$\trans A$ i~$A$ sovpadayut.
\endremark\cyrm

\proclaim{\cyrc\indent Lemma 7}\cyri
Pust\char"7E{} $\lambda\in\QQ$\rom, $b=\den\lambda\in\NN$ i $r\in\NN$.
Togda naimen\char"7Exi\char"1A{}
obwi\char"1A{} znamenatel\char"7E~$\psi_k$\rom, $k\in\NN$\rom, qisel
$$
\frac{\Delta_n^{(j)}(\lambda)}{j!},
\qquad j=0,1,\dots,r-1, \quad n=0,1,\dots,k,
$$
delit
$$
b^kd_k^{r-1}\prod_{p\mid b}p^{\tau_p(k)},
$$
gde $d_k$~-- naimen\char"7Exee obwee kratnoe qisel $1,2,\dots,k$,
a $\tau_p(k)$~-- stepen\char"7E{} vkhozhdeniya prostogo qisla~$p$ v~$k!$
\cyrm(sm.~\thetag{4}).
\endproclaim\cyrm

\demo{\cyrc\indent Dokazatel\char"7Estvo}\cyrm
Polagaya v teoreme raboty~\cite{13} \
$L=H=k$, $Q=b$, $x=-b\lambda$, $M=r-1$, $\Lambda=1$,
poluqaem trebuemoe.
\enddemo\cyrm

Iz lemmy~7 neposredstvenno sleduet

\proclaim{\cyrc\indent Lemma 8}\cyri
Pust\char"7E{} $b=\den(\lambda_1,\dots,\lambda_s)$~--
naimen\char"7Exi\char"1A{} obwi\char"1A{} znamenatel\char"7E{} qisel $\lambda_1,\dots,\lambda_s\in\QQ$\rom;
$r_1,\dots,r_s\in\NN$\rom, $r=\max_l\{r_l\}$.
Togda naimen\char"7Exi\char"1A{} obwi\char"1A{} znamenatel\char"7E{} $\psi_k$\rom, $k\in\NN$\rom,
qisel
$$
\frac{\Delta_n^{(j)}(\lambda_l)}{j!},
\qquad j=0,1,\dots,r_l-1, \quad l=1,\dots,s, \quad n=0,1,\dots,k,
$$
delit
$$
b^kd_k^{r-1}\prod_{p\mid b}p^{\tau_p(k)},
$$
gde $d_k$~-- naimen\char"7Exee obwee kratnoe qisel $1,2,\dots,k$\rom,
a $\tau_p(k)$~-- stepen\char"7E{} vkhozhdeniya prostogo qisla~$p$ v~$k!$\,.
\endproclaim\cyrm

Ob\char"7Fedinyaya rezul\char"7Etaty lemm~5,~6 i~8, poluqaem sleduyuwee
utverzhdenie.

\proclaim{\cyrc\indent Lemma 9}\cyri
Pust\char"7E~\thetag{17}~-- minimal\char"7Eny\char"1A{} mnogoqlen
racional\char"7Eno\char"1A{} matri\-cy~$A$\rom;
$r=\max_l\{r_l\}$\rom; $b=\den A$\rom;
$t_1,t_2$~-- naimen\char"7Exie obwie znamenateli e1lementov
matric~$T,T^{-1}$ sootvet\hbox{s}tvenno\rom,
gde $T$~-- matrica perekhoda ot~$A$ k ee zhordanovo\char"1A{} normal\char"7Eno\char"1A{} forme.
Togda naimen\char"7Exi\char"1A{} obwi\char"1A{} znamenatel\char"7E{} $\psi_k$\rom, $k\in\NN$\rom,
e1lementov matric
$$
\Delta_n(A), \qquad n=0,1,\dots,k,
$$
delit
$$
t_1t_2b^kd_k^{r-1}\prod_{p\mid b}p^{\tau_p(k)},
$$
gde $d_k$~-- naimen\char"7Exee obwee kratnoe qisel $1,2,\dots,k$\rom,
a $\tau_p(k)$~-- stepen\char"7E{} vkhozhdeniya prostogo qisla~$p$ v~$k!$\,.
\endproclaim\cyrm

S uqetom predel\char"7Enykh sootnoxeni\char"1A{}
$$
\limsup_{k\to\infty}d_k^{1/k}=e,
\qquad
\lim_{k\to\infty}\biggl(\prod_{p\mid b}p^{\tau_p(k)}\biggr)^{1/k}
=e^{\chi(b)}
\tag18
$$
i zameqaniya~2 k lemme~6 poluqaem sleduyuwee okonqatel\char"7Enoe utverzhdenie.

\proclaim{\cyrc\indent Teorema 3}\cyri
Pust\char"7E{} mnogoqlen
$$
P(\lambda)
=(\lambda-\lambda_1)^{r_1}
(\lambda-\lambda_2)^{r_2}\dotsb
(\lambda-\lambda_s)^{r_s},
\qquad \lambda_1,\dots,\lambda_s\in\QQ,
\quad \lambda_j\ne\lambda_l, \; j\ne l,
\tag19
$$
annuliruet matricu~$A$ \rom(naprimer\rom, $P(\lambda)$~-- minimal\char"7Eny\char"1A{}
ili kharakteristiqeski\char"1A{} mnogoqlen\rom)\rom;
$b$~-- naimen\char"7Exi\char"1A{} obwi\char"1A{} znamenatel\char"7E{} qisel $\lambda_1,\dots,\lambda_s$
\rom(znamenatel\char"7E{} matricy~$A$\rom);
$r=\max_l\{r_l\}$~-- maksimal\char"7Enaya kratnost\char"7E{} korne\char"1A{} mnogoqlena~\thetag{19}.
Togda matrica~$A$ udovletvoryaet usloviyu sokraweniya faktorialov
s postoyanno\char"1A{} $be^{\chi(b)+r-1}$\rom, gde funkciya $\chi(\,\cdot\,)$
opredelyaet\hbox{s}ya formulo\char"1A~\thetag{5}.
\endproclaim\cyrm

\subhead\cyrb
5. Sokrawenie faktorialov dlya sistem differencial\char"7Enykh
uravneni\char"1A{} fuksovskogo tipa
\endsubhead\cyrm
Vernemsya k voprosu, obsuzhdavxemusya v~p.~2.

Rassmotrim sistemu line\char"1Anykh differencial\char"7Enykh
uravneni\char"1A~\thetag{8} fuk\-sovskogo tipa i sootvet\hbox{s}tvuyuwie
sistemy~\thetag{10} dlya proizvodnykh poryadka $n$, $n=1,2,\dots$;
mnogoqlen $T(z)=(z-\gamma_1)\dotsb(z-\gamma_s)$
vyberem v sootvet\hbox{s}tvii s~\thetag{9}.
Dopolnim matricy neodnorodnykh sistem do kvadratnykh nulevymi strokami:
$$
\gather
Q(z)=\pmatrix
0 & 0 & \hdots & 0 \\
Q_{10}(z) & Q_{11}(z) & \hdots & Q_{1m}(z) \\
Q_{20}(z) & Q_{21}(z) & \hdots & Q_{2m}(z) \\
\hdotsfor4 \\
Q_{m0}(z) & Q_{m1}(z) & \hdots & Q_{mm}(z)
\endpmatrix,
\tag20
\\
Q^{[n]}(z)=\pmatrix
0 & 0 & \hdots & 0 \\
Q_{10}^{[n]}(z) & Q_{11}^{[n]}(z) & \hdots & Q_{1m}^{[n]}(z) \\
Q_{20}^{[n]}(z) & Q_{21}^{[n]}(z) & \hdots & Q_{2m}^{[n]}(z) \\
\hdotsfor4 \\
Q_{m0}^{[n]}(z) & Q_{m1}^{[n]}(z) & \hdots & Q_{mm}^{[n]}(z)
\endpmatrix,
\qquad n=1,2,\dots,
\endgather
$$
i perepixem rekurrentnye sootnoxeniya~\thetag{11}
v matriqnom vide:
$$
Q^{[n]}(z)=\frac d{dz}Q^{[n-1]}(z)+Q^{[n-1]}(z)Q(z),
\qquad n=1,2,\dots,
$$
otkuda
$$
\align
\trans Q^{[n]}(z)
&=\frac d{dz}\trans Q^{[n-1]}(z)+\trans Q(z)\trans Q^{[n-1]}(z)
=\biggl(\frac d{dz}+\trans Q(z)\biggr)\trans Q^{[n-1]}(z)
\\
&=\biggl(\frac d{dz}+\trans Q(z)\biggr)^nE,
\qquad n=1,2,\dots,
\tag21
\endalign
$$
gde $E$~-- ediniqnaya matrica razmera~$m+1$.

Ispol\char"7Ezuya razlozhenie~\thetag{12},
rassmotrim differencial\char"7Eny\char"1A{}
operator
$$
D=\frac d{dz}+\frac1{z-\gamma_1}\trans A_1+\dots
+\frac1{z-\gamma_s}\trans A_s,
\tag22
$$
gde $A_1,\dots,A_s$~-- racional\char"7Enye matricy.
Opredelim uslovie sokraweniya faktorialov dlya operatora~\thetag{22}\hlop,\
zamenyaya kol\char"7Eco $\ZZ_{\KK}[z]$ na $(\ZZ_{\KK}[z])^{(m+1)\times(m+1)}$
v opredelenii~1$'$.

\proclaim{\cyrc\indent Lemma 10}\cyri
Esli operator~\thetag{22} udovletvoryaet usloviyu
sokraweniya faktorialov s postoyanno\char"1A~$\Psi$\rom,
to sistema differencial\char"7Enykh uravneni\char"1A~\thetag{8}
tak\-zhe udovletvoryaet usloviyu sokraweniya faktorialov
s postoyanno\char"1A~$\Psi$.
\endproclaim\cyrm

\demo{\cyrc\indent Dokazatel\char"7Estvo}\cyrm
Soglasno sootnoxeniyam~\thetag{21} vypolneno
$$
\trans Q^{[n]}(z)=D^nE,
\qquad n=1,2,\dots\,.
\tag23
$$
Esli $\{\psi_k\}_{k\in\NN}$~-- posledovatel\char"7Enost\char"7E{}
iz opredeleniya sokraweniya faktorialov dlya operatora~$D$,
to iz~\thetag{23} sleduet, qto matricy
$$
\psi_k\frac{T^n(z)\trans Q^{[n]}(z)}{n!},
\qquad n=0,1,\dots,k, \quad k\in\NN,
$$
imeyut celoqislennye e1lementy.
Ot\hbox{s}yuda poluqaem utverzhdenie lemmy.
\enddemo\cyrm

K sozhaleniyu, lemma~4 neprimenima k
differencial\char"7Enomu operatoru \thetag{22}
dlya {\cyri proizvol\char"7Enykh\/} matric $A_1,\dots,A_s$;
tozhdestva
$$
\align
\frac{T^n(z)D^n}{n!}
&=\sum\Sb n_0,n_1,\dots,n_s\ge0\\n_0+n_1+\dots+n_s=n\endSb
(z-\gamma_1)^{n-n_1}\dotsb(z-\gamma_s)^{n-n_s}
\\ &\qquad\qquad\times
\frac{\<\trans A_1\>_{n_1}}{n_1!}\dotsb\frac{\<\trans A_s\>_{n_s}}{n_s!}
\cdot\frac1{n_0!}\frac{d^{n_0}}{dz^{n_0}},
\qquad n\in\NN,
\tag24
\endalign
$$
spravedlivy tol\char"7Eko v sluqae {\cyri poparno kommutiruyuwikh\/}
matric $A_1,\dots,A_s$.

Uqityvaya~\thetag{24}, lemmu~9 i ocenki~\thetag{18},
poluqaem sleduyuwee utverzhdenie.

\proclaim{\cyrc\indent Teorema 4}\cyri
Pust\char"7E{} $s=1$ ili matricy $A_1,\dots,A_s$
poparno kommutiruyut v sluqae $s\ge2$\rom;
$\lambda_1,\dots,\lambda_p\in\QQ$~-- sobstvennye znaqeniya
racional\char"7Enykh mat\-ric~$A_1,\dots,A_s$\rom;
$b=\den(\lambda_1,\dots,\lambda_p)$\rom;
$r_{jl}\ge0$\rom, $j=1,\dots,p$\rom, $l=1,\dots,s$\rom,~--
kratnost\char"7E{} sobstvennogo znaqeniya~$\lambda_j$
v minimal\char"7Enom \rom(kharakteristiqeskom\rom)
mnogoqlene matricy~$A_l$\rom;
$r=\max_{j,l}\{r_{jl}\}$~--
maksimal\char"7Enaya kratnost\char"7E{} sobstvennykh znaqeni\char"1A.
Togda operator~\thetag{22} udovletvoryaet usloviyu sokraweniya
faktorialov s postoyanno\char"1A{} $be^{\chi(b)+r-1}$\rom,
gde funkciya $\chi(\,\cdot\,)$ opredelena formulo\char"1A~\thetag{5}.
\endproclaim\cyrm

Soglasno lemme~10 iz teoremy~4 vytekaet sleduyuwee utverzhdenie.

\proclaim{\cyrc\indent Teorema 5}\cyri
Pust\char"7E{} matrica~\thetag{20}
sistemy line\char"1Anykh differencial\char"7Enykh
uravneni\char"1A~\thetag{8} imeet vid~\thetag{12}\rom,
gde $A_1,\dots,A_s$~-- racional\char"7Enye poparno kommutiruyuwie matricy.
Oboznaqim qerez~$b$ naimen\char"7Exi\char"1A{}
obwi\char"1A{} znamenatel\char"7E{} sobstvennykh znaqeni\char"1A{}
matric $A_1,\dots,A_s$\rom, qerez~$r$
maksimal\char"7Enuyu kratnost\char"7E{}
e1tikh sobstvennykh znaqeni\char"1A{} v minimal\char"7Enykh
\rom(kharakteristiqeskikh\rom) mnogoqlenakh
dannykh racional\char"7Enykh matric.
Togda sistema~\thetag{8} udovletvoryaet usloviyu sokraweniya faktorialov
s postoyanno\char"1A{} $be^{\chi(b)+r-1}$.
\endproclaim\cyrm

\remark{\cyrc\indent Zameqanie}\cyrm
V sluqae racional\char"7Enykh, no ne kommutiruyuwikh poparno matric
$A_1,\dots,A_s$ kaqestvennoe rexenie voprosa o sokrawenii
faktorialov dlya differencial\char"7Enogo operatora~\thetag{22}
ostaet\hbox{s}ya otkrytym. Qtoby poluqit\char"7E{} obobwenie
tozhdestv~\thetag{16} dlya operatora
$$
D=\frac d{dz}+\frac{A_1}{z-\gamma_1}+\dots+\frac{A_s}{z-\gamma_s}
\tag25
$$
opredelim matricy
$\<A_1,\dots,A_s\>_{\bn}$, $\bn=(n_1,\dots,n_s)$,
po indukcii, polagaya
$$
\align
&
\<A_1,A_2,\dots,A_s\>_{n_1,n_2,\dots,n_s}
\\ &\qquad
=\left\{\aligned
& 0, \quad \text{esli $\bn\notin(\ZZ_+)^s$}, \\
& E, \quad \text{esli $\bn=\0$}, \\
&
(A_1-n_1+1)\<A_1,A_2,\dots,A_s\>_{n_1-1,n_2,\dots,n_s}
\\ &\quad
+(A_2-n_2+1)\<A_1,A_2,\dots,A_s\>_{n_1,n_2-1,\dots,n_s}
+\dotsb
\\ &\quad
+(A_s-n_s+1)\<A_1,A_2,\dots,A_s\>_{n_1,n_2,\dots,n_s-1},
\quad \text{esli $\bn\in(\ZZ_+)^s$}.
\endaligned\right.
\endalign
$$
Esli matricy $A_1,\dots,A_s$ poparno kommutiruyut, to
$$
\gathered
\<A_1,\dots,A_s\>_{n_1,\dots,n_s}
=\frac{(n_1+\dots+n_s)!}{n_1!\dotsb n_s!}
\cdot\<A_1\>_{n_1}\dotsb\<A_s\>_{n_s},
\\
\bn=(n_1,\dots,n_s)\in(\ZZ_+)^s.
\endgathered
\tag26
$$

Indukcie\char"1A{} po~$k$ neslozhno dokazat\char"7E{} tozhdestva
$$
\<A_1+\dots+A_s\>_k
=\sum\Sb\bn\in(\ZZ_+)^s\\|\bn|=k\endSb
\<A_1,\dots,A_s\>_{\bn},
\qquad k=0,1,2,\dots,
$$
gde $|\bn|=n_1+\dots+n_s$, i
$$
D^k=\sum_{l=0}^k\binom kl
\Biggl(\sum\Sb\bn\in(\ZZ_+)^s\\|\bn|=k-l\endSb
\frac{\<A_1,\dots,A_s\>_{\bn}}
{(z-\gamma_1)^{n_1}\dotsb(z-\gamma_s)^{n_s}}\Biggr)
\frac{d^l}{dz^l},
\qquad k=0,1,2,\dots,
$$
dlya differencial\char"7Enogo operatora~\thetag{25}.
V qastnosti,
$$
D^kE
=\sum\Sb\bn\in(\ZZ_+)^s\\|\bn|=k\endSb
\frac{\<A_1,\dots,A_s\>_{\bn}}
{(z-\gamma_1)^{n_1}\dotsb(z-\gamma_s)^{n_s}},
\qquad k=0,1,2,\dots\,.
\tag27
$$
\endremark\cyrm

\subhead\cyrb
6. Prilozhenie k obobwennym gipergeometriqeskim funkciyam
\endsubhead\cyrm
Ob\-obwennaya gipergeometriqeskaya funkciya
$$
\gather
f(z)=F\biggl(\gathered
\alpha_1,\dots,\alpha_m \\ \beta_1+1,\dots,\beta_m+1
\endgathered\Bigm|z\biggr)
=\sum_{n=0}^\infty\frac{\<-\alpha_1\>_n\dotsb\<-\alpha_m\>_n}
{\<-\beta_1-1\>_n\dotsb\<-\beta_m-1\>_n}z^n,
\tag28
\\
\beta_1,\dots,\beta_m\notin\{-1,-2,\dots\},
\endgather
$$
udovletvoryaet line\char"1Anomu differencial\char"7Enomu uravneniyu
\cite{14, gl.~5, \S\,1, lemma~1 pri $t=1$}
$$
\bigl((\delta+\beta_1)\dotsb(\delta+\beta_m)
-z(\delta+\alpha_1)\dotsb(\delta+\alpha_m)\bigr)y
=\beta_1\dotsb\beta_m,
\qquad \delta=z\frac d{dz},
$$
poryadka~$m$.

Dlya funkci\char"1A{}
$$
f_1(z)=f(z), \quad f_2(z)=\delta f_1(z), \quad \dots,
\quad f_m(z)=\delta f_{m-1}(z)
\tag29
$$
poluqaem sistemu line\char"1Anykh differencial\char"7Enykh uravneni\char"1A{}
$$
\gathered
\frac d{dz}y_l=\frac1zy_{l+1}, \qquad l=1,\dots,m-1,
\\
\aligned
\frac d{dz}y_m
&=\frac{\sigma_1(\bbeta)-z\sigma_1(\balpha)}{z(z-1)}y_m
+\frac{\sigma_2(\bbeta)-z\sigma_2(\balpha)}{z(z-1)}y_{m-1}
+\dotsb
\\ \noalign{\vskip-1pt} &\qquad
+\frac{\sigma_m(\bbeta)-z\sigma_m(\balpha)}{z(z-1)}y_1
-\frac{\sigma_m(\bbeta)}{z(z-1)},
\endaligned
\endgathered
\tag30
$$
gde $\sigma_l(\,\cdot\,)$, $l=1,\dots,m$,~-- simmetriqeskie
mnogoqleny Vieta stepeni~$l$, t.e\.
$$
\gathered
(z+\alpha_1)(z+\alpha_2)\dotsb(z+\alpha_m)
=z^m+\sigma_1(\balpha)z^{m-1}+\dots
+\sigma_{m-1}(\balpha)z+\sigma_m(\balpha),
\\
(z+\beta_1)(z+\beta_2)\dotsb(z+\beta_m)
=z^m+\sigma_1(\bbeta)z^{m-1}+\dots
+\sigma_{m-1}(\bbeta)z+\sigma_m(\bbeta).
\endgathered
\tag31
$$
S uqetom
$$
\gather
\frac{\sigma_l(\bbeta)-z\sigma_l(\balpha)}{z(z-1)}
=\frac{\sigma_l(\bbeta)-\sigma_l(\balpha)}{z-1}
-\frac{\sigma_l(\bbeta)}z, \qquad l=1,\dots,m,
\\
\frac1{z(z-1)}=\frac1{z-1}-\frac1z
\endgather
$$
perepixem sistemu differencial\char"7Enykh uravneni\char"1A~\thetag{30}
v matriqnom vide:
$$
\align
&
\frac d{dz}\pmatrix y_1 \\ y_2 \\ \hdots \\ y_{m-1} \\ y_m \endpmatrix
=\frac1z\pmatrix 0 \\ 0 \\ \hdots \\ 0 \\ \sigma_m(\bbeta) \endpmatrix
+\frac1{z-1}\pmatrix 0 \\ 0 \\ \hdots \\ 0 \\ -\sigma_m(\bbeta) \endpmatrix
\\ &\;\;
+\frac1z\pmatrix
0 & 1 & 0 & \hdots & 0 \\
0 & 0 & 1 & \hdots & 0 \\
\hdotsfor5 \\
0 & 0 & \hdots & 0 & 1 \\
-\sigma_m(\bbeta) & -\sigma_{m-1}(\bbeta) & \hdots &
-\sigma_2(\bbeta) & -\sigma_1(\bbeta)
\endpmatrix
\pmatrix y_1 \\ y_2 \\ \hdots \\ y_{m-1} \\ y_m \endpmatrix
\\ &\;\;
+\frac1{z-1}\pmatrix
0 & \hdots & 0 & 0 \\
0 & \hdots & 0 & 0 \\
\hdotsfor4 \\
0 & \hdots & 0 & 0 \\
\sigma_m(\bbeta)-\sigma_m(\balpha) & \hdots &
\sigma_2(\bbeta)-\sigma_2(\balpha) &
\sigma_1(\bbeta)-\sigma_1(\balpha)
\endpmatrix
\pmatrix y_1 \\ y_2 \\ \hdots \\ y_{m-1} \\ y_m \endpmatrix.
\\ \noalign{\vskip-3pt}
\tag32
\endalign
$$

Teorema~5 k sisteme~\thetag{32} neprimenima, tak kak
otveqayuwie regul\hbox{ya}r\-nym osobennostyam $z=0$ i $z=1$ matricy
ne kommutiruyut. No my i ne stavim celi
{\cyri sokratit\char"7E{} faktorialy\/} dlya sistemy~\thetag{32}.
Toqnoe vyqislenie postoyanno\char"1A, s kotoro\char"1A{} proiskhodit
sokrawenie faktorialov dlya gipergeometriqeskogo differencial\char"7Enogo
uravneniya, posle raboty~\cite{9} stalo neaktual\char"7Enym:
ispol\char"7Ezovanie priblizheni\char"1A{} Pade vtorogo roda
(vmesto pervogo) v kaqestve priblizhayuwikh funkcional\char"7Enykh form
daet vozmozhnost\char"7E{} sokrawat\char"7E{} faktorialy bez dopolnitel\char"7Enykh
usili\char"1A{} (sm.~\cite{9} i~\cite{15}).
Naxa cel\char"7E~-- posqitat\char"7E{} postoyannuyu, s kotoro\char"1A{} proiskhodit
sokrawenie faktorialov dlya sistemy line\char"1Anykh differencial\char"7Enykh
uravneni\char"1A, sopryazhenno\char"1A{} k odnorodno\char"1A{} qasti sistemy~\thetag{30};
e1to dast vozmozhnost\char"7E{} vospo\-l\char"7Ezovat\char"7Esya
osnovno\char"1A{} teoremo\char"1A{} iz~\cite{16}
i poluqit\char"7E{} ocenki mery irraciona\-l\char"7Enosti
znaqeni\char"1A{} gipergeometriqesko\char"1A{} funkcii~\thetag{28}
(bez posledovatel\char"7E\-nykh proizvodnykh)
v racional\char"7Enykh toqkakh.

Sopryazhennaya sistema dlya obwego sluqaya~\thetag{8}
poluqaet\hbox{s}ya transponiro\-vaniem i smeno\char"1A{} znaka matricy
$\bigl(Q_{lj}(z)\bigr)_{l,j=1,\dots,m}$; v naxem sluqae
soprya\-zhennaya sistema imeet vid
$$
\frac d{dz}\pmatrix y_1 \\ y_2 \\ \hdots \\ y_{m-1} \\ y_m \endpmatrix
=\biggl(\frac1zA_1+\frac1{z-1}A_2\biggr)
\pmatrix y_1 \\ y_2 \\ \hdots \\ y_{m-1} \\ y_m \endpmatrix,
\tag33
$$
gde
$$
\allowdisplaybreaks
\align
A_1&=\pmatrix
0 & 0 & \hdots & 0 & \sigma_m(\bbeta) \\
-1 & 0 & \hdots & 0 & \sigma_{m-1}(\bbeta) \\
0 & -1 & \hdots & 0 & \sigma_{m-2}(\bbeta) \\
\hdotsfor5 \\
0 & \hdots & -1 & 0 & \sigma_2(\bbeta) \\
0 & \hdots & 0 & -1 & \sigma_1(\bbeta)
\endpmatrix,
\\
A_2&=\pmatrix
0 & \hdots & 0 & \sigma_m(\balpha)-\sigma_m(\bbeta) \\
0 & \hdots & 0 & \sigma_{m-1}(\balpha)-\sigma_{m-1}(\bbeta) \\
\hdotsfor4 \\
0 & \hdots & 0 & \sigma_2(\balpha)-\sigma_2(\bbeta) \\
0 & \hdots & 0 & \sigma_1(\balpha)-\sigma_1(\bbeta)
\endpmatrix.
\endalign
$$

Pervye $m-1$ stolbcov matricy~$A_2$ nulevye, t.e\. rang e1to\char"1A{} matricy
raven~$m-1$. Sledovatel\char"7Eno, $m-1$ sobstvennykh znaqeni\char"1A{} $\lambda=0$
matricy~$A_2$ vkhodyat s kratnost\char"7Eyu~1 v minimal\char"7Eny\char"1A{} mnogoqlen.
Ostavxeesya sobstvennoe znaqenie sovpadaet
so sledom matricy~$A_2$ i ravno
$\gamma=\sigma_1(\balpha)-\sigma_1(\bbeta)
=\allowmathbreak\alpha_1+\dots+\alpha_m-\beta_1-\dotsb-\beta_m$.

Matrica~$A_1$ s toqnost\char"7Eyu do znaka yavlyaet\hbox{s}ya
{\cyri kletko\char"1A{} Frobeniusa\/}
(sm.~\cite{12, gl.~6, \S\,6}); ee kharakteristiqeski\char"1A{} mnogoqlen raven
$$
\align
\det(A_1-\lambda E)
&=(-1)^m(\lambda^m-\sigma_1\lambda^{m-1}
+\sigma_2\lambda^{m-2}+\dots+(-1)^m\sigma_m)
\\
&=(-1)^m(\lambda-\beta_1)(\lambda-\beta_2)\dotsb(\lambda-\beta_m).
\endalign
$$
Poe1tomu sobstvennye znaqeniya matricy~$A_1$ ravny $\beta_1,\dots,\beta_m$.

\proclaim{\cyrc\indent Lemma 11}\cyri
Pust\char"7E{} parametry $\beta_1,\dots,\beta_m$
poparno razliqny i $\gamma\ne0$\rom,
$b$~-- naimen\char"7Exi\char"1A{} obwi\char"1A{}
znamenatel\char"7E{} qisel $\gamma,\beta_1,\dots,\beta_m$.
Togda sistema~\thetag{33} udov\-letvoryaet usloviyu sokraweniya
faktorialov s postoyanno\char"1A{} $be^{\chi(b)+2}$\rom,
gde funk\-ciya $\chi(\,\cdot\,)$
opredelyaet\hbox{s}ya formulo\char"1A~\thetag{5}.
\endproclaim\cyrm

\demo{\cyrc\indent Dokazatel\char"7Estvo}\cyrm
Sobstvennomu znaqeniyu~$\beta_j$ matricy~$A_1$ otveqaet
so\-bstvenny\char"1A{} vektor
$$
\pmatrix t_{1j} \\ t_{2j} \\ \hdots \\ t_{mj} \endpmatrix
=\pmatrix
\sigma_{m-1}(\beta_1,\dots,\beta_{j-1},\beta_{j+1},\dots,\beta_m) \\
\sigma_{m-2}(\beta_1,\dots,\beta_{j-1},\beta_{j+1},\dots,\beta_m) \\
\hdotsfor1 \\
\sigma_1(\beta_1,\dots,\beta_{j-1},\beta_{j+1},\dots,\beta_m) \\
\sigma_0(\beta_1,\dots,\beta_{j-1},\beta_{j+1},\dots,\beta_m)
\endpmatrix,
\qquad j=1,\dots,m,
$$
gde $\sigma_l(\,\cdot\,)$~-- koe1fficient pri~$z^{m-1-l}$ v mnogoqlene
$$
(z+\beta_1)\dotsb(z+\beta_{j-1})(z+\beta_{j+1})\dotsb(z+\beta_m).
$$
Dlya matricy perekhoda $T=(t_{lj})_{l,j=1,\dots,m}$
polozhim $\wt T=T^{-1}=(\tilde t_{lj})_{l,j=1,\dots,m}$.
Togda
$$
\tilde t_{lj}
=(-1)^{m+j}\beta_l^{j-1}
\cdot\prod\Sb k=1\\k\ne l\endSb^m\frac1{\beta_l-\beta_k},
\qquad l,j=1,\dots,m,
$$
matrica $\wt A_1=T^{-1}A_1T$ diagonal\char"7Enaya,
$\wt A_1=\diag(\beta_1,\dots,\beta_m)$,
a matrica $\wt A_2=\allowmathbreak T^{-1}A_2T$ imeet vid
$$
\wt A_2=\pmatrix a_1 \\ \hdots \\ a_m \endpmatrix
(1 \; \hdots \; 1),
\qquad
a_j
=-(\beta_j-\alpha_j)
\cdot\prod\Sb k=1\\k\ne j\endSb^m\frac{\beta_j-\alpha_k}{\beta_j-\beta_k},
\quad j=1,\dots,m.
$$

Zapixem matricy $Q^{[n]}(z)$, $n=0,1,2,\dots$, iz opredeleniya~3
dlya sistemy differencial\char"7Enykh uravneni\char"1A~\thetag{33}.
Soglasno~\thetag{23} i~\thetag{27} imeem
$$
\align
\trans Q^{[n]}(z)
&=\biggl(\frac d{dz}+\frac1z\trans A_1+\frac1{z-1}\trans A_2\biggr)^nE
\\
&=\sum\Sb n_1,n_2\ge 0\\n_1+n_2=n\endSb
\frac{\<\trans A_1,\trans A_2\>_{n_1,n_2}}
{z^{n_1}(z-1)^{n_2}},
\qquad n=0,1,2,\dots,
\tag34
\endalign
$$
gde
$$
\<\trans A_1,\trans A_2\>_{n_1,n_2}
=\left\{\aligned
& 0, \quad \text{esli $n_1<0$ ili $n_2<0$}, \\
& E, \quad \text{esli $n_1=n_2=0$}, \\
& (\trans A_1-n_1+1)\<\trans A_1,\trans A_2\>_{n_1-1,n_2}
\\ &\quad
+(\trans A_2-n_2+1)\<\trans A_1,\trans A_2\>_{n_1,n_2-1}
\quad \text{v protivnom sluqae}
\endaligned\right.
$$
(sm\. zameqanie k teoreme~5).
Nekommutativnost\char"7E{} matric $\trans A_1,\trans A_2$
v tozhde\-stvakh~\thetag{34} ne spasayut dazhe sootnoxeniya
$$
\gather
\frac1{z^{n_1+1}(1-z)^{n_2+1}}
=\sum_{k=0}^{n_1}\binom{n_1+n_2-k}{n_2}\frac1{z^{k+1}}
+\sum_{k=0}^{n_2}\binom{n_1+n_2-k}{n_1}\frac1{(1-z)^{k+1}},
\\
n_1,n_2=0,1,2,\dots\,.
\tag35
\endgather
$$

Polagaya
$$
B_1=\trans{\wt A_1}=\diag(\beta_1,\dots,\beta_m),
\qquad
B_2=\trans{\wt A_2}
=\pmatrix 1 \\ \hdots \\ 1 \endpmatrix
(a_1 \; \hdots \; a_m),
$$
iz~\thetag{34} poluqaem
$$
\trans{\bigl(T^{-1}Q^{[n]}(z)T\bigr)}
=\sum\Sb n_1,n_2\ge 0\\n_1+n_2=n\endSb
\frac{\<B_1,B_2\>_{n_1,n_2}}
{z^{n_1}(z-1)^{n_2}},
\qquad n=0,1,2,\dots\,.
\tag36
$$

Zafiksiruem teper\char"7E{} paru celykh neotricatel\char"7Enykh $n_1,n_2$.
V sluqae $n_2=0$ znamenatel\char"7E{} matricy
$$
\frac{\<B_1,B_2\>_{n_1,n_2}}{(n_1+n_2)!}
=\frac{\<B_1\>_{n_1}}{n_1!}
$$
soglasno lemme~9 delit $b_0^k\prod_{p\mid b_0}p^{\tau_p(k)}$
dlya lyubogo $k\ge n_1=n_1+n_2$;
zdes\char"7E{} $b_0=\allowmathbreak\den(\beta_1,\dots,\beta_m)$.

Pust\char"7E, dalee, $n_2>0$.
Soglasno rekurrentnym sootnoxeniyam dlya mat\-ric $\<B_1,B_2\>_{n_1,n_2}$
(i formule~\thetag{26} v sluqae, kogda matricy kommutiruyut)
zaklyuqaem, qto matrica $\<B_1,B_2\>_{n_1,n_2}$
predstavlyaet sobo\char"1A{}
summu $N=\mathbreak(n_1+n_2)!/(n_1!\,n_2!)$
slagaemykh, kazhdoe iz kotorykh s toqnost\char"7Eyu
do poryadka mnozhitele\char"1A{}
sovpadaet s~$\<B_1\>_{n_1}\<B_2\>_{n_2}$:
$$
\<B_1,B_2\>_{n_1,n_2}=\sum_{r=1}^NB^{(r)}.
\tag37
$$
Srazu otmetim, qto esli khotya by odna iz kvadratnykh matric~$X_1$ i~$X_2$
diagonal\char"7Ena, to glavnye diagonali matric $X_1X_2$ i~$X_2X_1$ sovpadayut.
Posko\-l\char"7Eku matricy $B_1-lE$, $l=0,1,\dots,n_1-1$, diagonal\char"7Eny,
glavnaya diagonal\char"7E{} kazhdogo slagaemogo v~\thetag{37} sovpadaet
s glavno\char"1A{} diagonal\char"7Eyu matricy $\<B_1\>_{n_1}\<B_2\>_{n_2}$.
Matrica~$B_2$ udovletvoryaet sootnoxeniyu
$$
B_2^2=(a_1+\dots+a_m)B_2=\Tr B_2\cdot B_2=\gamma B_2,
$$
otkuda $\<B_2\>_{n_2}=\<\gamma-1\>_{n_2-1}B_2$ i
$$
\<B_1\>_{n_1}\<B_2\>_{n_2}
=\<\gamma-1\>_{n_2-1}\cdot\pmatrix
a_1\<\beta_1\>_{n_1} & a_2\<\beta_1\>_{n_1} & \hdots & a_m\<\beta_1\>_{n_1} \\
a_1\<\beta_2\>_{n_1} & a_2\<\beta_2\>_{n_1} & \hdots & a_m\<\beta_2\>_{n_1} \\
\hdotsfor4 \\
a_1\<\beta_m\>_{n_1} & a_2\<\beta_m\>_{n_1} & \hdots & a_m\<\beta_m\>_{n_1}
\endpmatrix.
\tag38
$$

Rassmotrim lyuboe iz slagaemykh $B=B^{(r)}$, $1\le r\le N$,
vkhodyawikh v~\thetag{37}. Mnozhitel\char"7E{} $B_1-jE$ v e1tom slagaemom poyavlyaet\hbox{s}ya
levee (t.e\. pozdnee) mnozhitelya $B_1-lE$ togda i tol\char"7Eko togda,
kogda $j>l$; to zhe samoe mozhno skazat\char"7E{} o poryadke poyavleniya v~$B$
mnozhitele\char"1A{} $B_2-jE$ i $B_2-lE$. Vydelyaya pervoe (i edinstvennoe) poyavlenie
mnozhitelya~$B_2$ v slagaemom~$B$, nakhodim
$$
\align
B
&=X\cdot B_2\cdot\<B_1\>_s
=X\cdot\pmatrix 1 \\ \hdots \\ 1 \endpmatrix
(a_1 \; \hdots \; a_m)
\cdot\diag\bigl(\<\beta_1\>_s,\dots,\<\beta_m\>_s\bigr)
\\
&=\pmatrix x_1 \\ x_2 \\ \hdots \\ x_m \endpmatrix
\bigl(a_1\<\beta_1\>_s \; a_2\<\beta_2\>_s \;
\hdots \; a_m\<\beta_m\>_s\bigr)
\\
&=\pmatrix
a_1\<\beta_1\>_sx_1 & a_2\<\beta_2\>_sx_1 & \hdots & a_m\<\beta_m\>_sx_1 \\
a_1\<\beta_1\>_sx_2 & a_2\<\beta_2\>_sx_2 & \hdots & a_m\<\beta_m\>_sx_2 \\
\hdotsfor4 \\
a_1\<\beta_1\>_sx_m & a_2\<\beta_2\>_sx_m & \hdots & a_m\<\beta_m\>_sx_m
\endpmatrix.
\tag39
\endalign
$$
Sravnenie glavnykh diagonale\char"1A{} matric~\thetag{38} i~\thetag{39} daet
$$
x_l=\<\gamma-1\>_{n_2-1}\cdot\<\beta_l-s\>_{n_1-s},
\qquad l=1,\dots,m.
$$
Takim obrazom, kazhdoe slagaemoe, vkhodyawee v~\thetag{37}, imeet vid
$$
B=\<\gamma-1\>_{n_2-1}\cdot\Bigl(
a_j\<\beta_j\>_s\<\beta_l-s\>_{n_1-s}
\Bigr)_{l,j=1,\dots,m}
$$
dlya nekotorogo~$s$, $0\le s<n_1$.
Polagaya $a=\den(a_1,\dots,a_m)$,
$b=\den(\gamma,\beta_1,\dots,\allowmathbreak\beta_m)$,
$g_k$~-- naimen\char"7Exee obwee kratnoe qisel
$$
\frac{k!}{k_0!\,k_1!\,k_2!},
\qquad k_0,k_1,k_2=0,1,2,\dots, \quad k_0+k_1+k_2=k,
\tag40
$$
i uqityvaya uslovie $\gamma\ne0$,
soglasno lemme~9 poluqaem, qto naimen\char"7Exi\char"1A{} obwi\char"1A{} znamenatel\char"7E{}
e1lementov matricy
$$
\frac{\gamma B}{(n_1+n_2)!}
=\frac{s!\,(n_1-s)!\,n_2!}{(n_1+n_2)!}
\cdot\frac{\gamma B}{s!\,(n_1-s)!\,n_2!}
$$
delit
$$
g_k\cdot ab^k\prod_{p\mid b}p^{\tau_p(k)}
\tag41
$$
dlya lyubogo $k\ge n_1+n_2$; vvidu proizvol\char"7Enosti vybora slagaemogo~$B$
v summe~\thetag{37} zaklyuqaem, qto naimen\char"7Exi\char"1A{} obwi\char"1A{} znamenatel\char"7E{}
e1lementov mat\-ricy
$$
\frac{\gamma\<B_1,B_2\>_{n_1,n_2}}{(n_1+n_2)!}
=\sum_{r=1}^N\frac{\gamma B^{(r)}}{(n_1+n_2)!}
$$
takzhe delit qislo~\thetag{41} dlya lyubogo $k\ge n_1+n_2$.

Stepen\char"7E{} vkhozhdeniya prostogo~$p$ v kazhdoe qislo~\thetag{40}
soglasno~\thetag{4} ravna
$$
\gather
\tau_p(k)-\tau_p(k_0)-\tau_p(k_1)-\tau_p(k_2)
=\sum_{m=1}^\infty\biggl(\biggl[\frac k{p^m}\biggr]
-\biggl[\frac{k_0}{p^m}\biggr]-\biggl[\frac{k_1}{p^m}\biggr]
-\biggl[\frac{k_2}{p^m}\biggr]\biggr),
\tag42
\\
k_0,k_1,k_2=0,1,2,\dots, \quad k_0+k_1+k_2=k.
\endgather
$$
Summirovanie v~\thetag{42} proiskhodit tol\char"7Eko po $m\le[\log k/\log p]$;
krome togo,
$$
[\xi_0+\xi_1+\xi_2]-[\xi_1]-[\xi_2]-[\xi_3]\le2,
\qquad \xi_0,\xi_1,\xi_2\in\Bbb R.
$$
Sledovatel\char"7Eno, dlya vsekh prostykh~$p\le k$ vypolneno
$$
\gather
\tau_p(k)-\tau_p(k_0)-\tau_p(k_1)-\tau_p(k_2)
\le2\biggl[\frac{\log k}{\log p}\biggr]
\le2\frac{\log k}{\log p},
\\
k_0,k_1,k_2=0,1,2,\dots, \quad k_0+k_1+k_2=k.
\endgather
$$
E1ti neravenstva dayut ocenku
$$
g_k\le\prod_{p\le k}p^{2\log k/\log p}=e^{2\pi(k)\log k},
$$
gde $\pi(k)$~-- koliqestvo prostykh, ne prevoskhodyawikh~$k$, otkuda
$$
\limsup_{k\to\infty}g_k^{1/k}\le e^2.
\tag43
$$
S uqetom poluqennogo predel\char"7Enogo sootnoxeniya,
tozhdestv~\thetag{36} i lemmy~5 my poluqaem, qto
sistema odnorodnykh differencial\char"7Enykh uravneni\char"1A~\thetag{33}
udovletvoryaet usloviyu sokraweniya faktorialov s postoyanno\char"1A{}
$be^{\chi(b)+2}$. Lemma dokazana.
\enddemo\cyrm

\remark{\cyrc\indent Zameqanie 1}\cyrm
Povtoryaya privedennye rassuzhdeniya v sluqae $\gamma=0$,
dlya lyubogo slagaemogo v~\thetag{37} poluqaem
$$
B=(-1)^{n_2-1}(n_2-1)!\cdot\Bigl(
a_j\<\beta_j\>_s\<\beta_l-s\>_{n_1-s}
\Bigr)_{l,j=1,\dots,m}
$$
s nekotorym~$s$, $0\le s<n_1$, otkuda naimen\char"7Exi\char"1A{}
obwi\char"1A{} znamenatel\char"7E{} e1lemen\-tov matricy
$$
\frac{\<B_1,B_2\>_{n_1,n_2}}{(n_1+n_2)!}
$$
dlya lyubogo $k\ge n_1+n_2$ delit
$$
d_kg_k\cdot ab^k\prod_{p\mid b}p^{\tau_p(k)},
$$
gde $d_k$~-- naimen\char"7Exee obwee kratnoe qisel $1,2,\dots,k$.
S uqetom predel\char"7Enykh sootnoxeni\char"1A~\thetag{18} i~\thetag{43}
sokrawenie faktorialov v e1tom sluqae proiskho\-dit s postoyanno\char"1A{}
$be^{\chi(b)+3}$.
\endremark\cyrm

\remark{\cyrc\indent Zameqanie 2}\cyrm
Dokazatel\char"7Estvo lemmy~11 godit\hbox{s}ya
dlya podsqeta (to\char"1A{} zhe samo\char"1A)
postoyanno\char"1A, s kotoro\char"1A{} proiskhodit sokrawenie faktorialov
u iskhodno\char"1A{} neodnorodno\char"1A{} sistemy~\thetag{32} v sluqae
poparno razliqnykh para\-metrov $\beta_1,\dots,\beta_m$.
\endremark\cyrm

Sleduyuwie utverzhdeniya otnosyat\hbox{s}ya k arifmetiqeskim
i algebraiqes\-kim svo\char"1Astvam funkci\char"1A~\thetag{29}.

\proclaim{\cyrc\indent Lemma 12}\cyri
Pust\char"7E{} $b_1,b_2$~-- znamenateli qisel
$\alpha,\beta\in\QQ\setminus\{-1,-2,\dots\}$\rom;
$b$~-- naimen\char"7Exee obwee kratnoe qisel~$b_1,b_2$.
Oboznaqim qerez $\phi_k\in\NN$\rom, $k\in\NN$\rom,
naimen\char"7Exi\char"1A{} obwi\char"1A{}
znamenatel\char"7E{} racional\char"7Enykh qisel
$$
\frac{\<-\alpha\>_n}{\<-\beta\>_n}, \qquad n=0,1,\dots,k.
$$
Togda
$$
\limsup_{k\to\infty}\phi_k^{1/k}
\le\Phi=e^{\rho(b_2)}\frac{b_1}b,
$$
gde
$$
\rho(b)=\frac b{\phi(b)}\sum\Sb 1\le n\le b\\(n,b)=1\endSb\frac1n,
\quad
\phi(b)=\sum\Sb 1\le n\le b\\(n,b)=1\endSb1,
\qquad b\in\NN.
\tag44
$$
\endproclaim\cyrm

\demo{\cyrc\indent Dokazatel\char"7Estvo}\cyrm
Naimen\char"7Exie obwie znamenateli qisel
$$
\<-\alpha\>_n, \quad n=0,1,\dots,k,
\qquad
\<-\beta\>_n, \quad n=0,1,\dots,k,
$$
ravny sootvet\hbox{s}tvenno $b_1^k$ i $b_2^k$.
Naimen\char"7Exi\char"1A{} obwi\char"1A{} znamenatel\char"7E{} qisel
$$
\frac1{\<-\beta\>_n}, \qquad n=0,1,\dots,k,
$$
raven naimen\char"7Exemu obwemu kratnomu~$d_k$ qisel
$$
-a+b_2(n-1), \quad n=1,\dots,k,
\qquad a=b_2\beta\in\ZZ.
$$
Soglasno~\cite{16, lemma~3.2} spravedliva ocenka
$$
\limsup_{k\to\infty}\frac{\log d_k}k
\le\rho(b_2),
$$
gde funkciya $\rho(\,\cdot\,)$ opredelyaet\hbox{s}ya ravenstvom~\thetag{44}.
E1to dokazyvaet lemmu.
\enddemo\cyrm

\proclaim{\cyrc\indent Lemma 13}\cyri
Pust\char"7E{} $q_1,q_2$~-- proizvedeniya znamenatele\char"1A{} qisel
$\alpha_1,\dots,\alpha_m$ i $\beta_1,\dots,\beta_m$ sootvet\hbox{s}tvenno\rom;
$b$~-- naimen\char"7Exee obwee kratnoe qisel~$q_1,q_2$\rom;
$b_1,\dots,b_m$~-- znamenateli qisel
$\beta_1,\dots,\beta_m$ sootvet\hbox{s}tvenno.
Togda funkciya \thetag{28} prinadlezhit klassu
$\GG(1,\Phi)$\rom, gde $\Phi=e^{\rho(b_1)+\dots+\rho(b_m)}q_1/b$.
\endproclaim\cyrm

\demo{\cyrc\indent Dokazatel\char"7Estvo}\cyrm
Ocenka veliqiny~$\Phi$ vytekaet iz lemmy~12.
Poskol\char"7E\-ku
$$
\lim_{n\to\infty}\biggl|\frac{\<-\lambda\>_n}{n!}\biggr|^{1/n}=1,
\qquad \lambda\notin\{-1,-2,\dots\},
$$
oblast\char"7E{} skhodimosti ryada~\thetag{28}~-- krug $|z|<1$.
Lemma dokazana.
\enddemo\cyrm

\proclaim{\cyrc\indent Lemma 14}\cyri
Esli $f(z)\in\GG(C,\Phi)$\rom,
to $\delta f(z)\in\GG(C,\Phi)$\rom, gde $\delta=z\frac d{dz}$.
\endproclaim\cyrm

\demo{\cyrc\indent Dokazatel\char"7Estvo}\cyrm
E1to utverzhdenie oqevidno. Esli ryad Te\char"1Alora
$$
f(z)=\sum_{n=0}^\infty f_nz^n
$$
funkcii $f(z)$ skhodit\hbox{s}ya v kruge $|z|<C$, to v e1tom zhe kruge
skhodit\hbox{s}ya ryad
$$
\delta f(z)=\sum_{n=1}^\infty nf_nz^n.
$$
Esli posledovatel\char"7Enost\char"7E{} natural\char"7Enykh
qisel $\{\phi_k\}_{k\in\NN}$ vybrana tak, qto
$$
\phi_kf_n\in\ZZ, \qquad n=0,1,\dots,k, \quad k\in\NN,
$$
to
$$
\phi_knf_n\in\ZZ, \qquad n=0,1,\dots,k, \quad k\in\NN.
$$
Takim obrazom, $\delta f(z)\in\GG(C,\Phi)$, qto i trebovalos\char"7E{}
dokazat\char"7E.
\enddemo\cyrm

\proclaim{\cyrc\indent Sledstvie}\cyri
Pust\char"7E{} $q_1,q_2$~-- proizvedeniya znamenatele\char"1A{} qisel
$\alpha_1,\dots,\alpha_m$ i $\beta_1,\dots,\beta_m$ sootvet\hbox{s}tvenno\rom;
$b$~-- naimen\char"7Exee obwee kratnoe qisel~$q_1,q_2$\rom;
$b_1,\dots,b_m$~-- znamenateli qisel $\beta_1,\dots,\beta_m$
sootvet\hbox{s}tvenno.
Togda sovokupnost\char"7E{} funkci\char"1A~\thetag{29} prinadlezhit klassu
$\GG(1,\Phi)$\rom, gde $\Phi=e^{\rho(b_1)+\dots+\rho(b_m)}q_1/b$.
\endproclaim\cyrm

\proclaim{\cyrc\indent Lemma 15}\cyri
Vronskian \rom(opredelitel\char"7E{} matricy
fundamental\char"7Eno\char"1A{} siste\-my rexeni\char"1A\rom)
$W(z)$ line\char"1Anogo odnorodnogo differencial\char"7Enogo
uravneniya
$$
\bigl((\delta+\beta_1)\dotsb(\delta+\beta_m)
-z(\delta+\alpha_1)\dotsb(\delta+\alpha_m)\bigr)y=0,
\qquad \delta=z\frac d{dz},
\tag45
$$
poryadka~$m$ udovletvoryaet differencial\char"7Enomu uravneniyu
$$
\bigl((\delta+\beta)-z(\delta+\alpha)\bigr)y=0,
\qquad \alpha=\sigma_1(\balpha), \quad \beta=\sigma_1(\bbeta),
\tag46
$$
i\rom, sledovatel\char"7Eno\rom, dlya racional\char"7Enykh~$\alpha,\beta$
yavlyaet\hbox{s}ya algebraiqesko\char"1A{} funkcie\char"1A\rom:
$$
W(z)=Cz^{-\beta}(1-z)^{\alpha-\beta}, \qquad C\in\CC.
\tag47
$$
\endproclaim\cyrm

\demo{\cyrc\indent Dokazatel\char"7Estvo}\cyrm
Vronskian~$W(z)$ differencial\char"7Enogo uravneniya~\thetag{45} sovpadaet
s vronskianom sistemy line\char"1Anykh differencial\char"7Enykh uravneni\char"1A{}
$$
\gathered
\frac d{dz}y_l=\frac1zy_{l+1}, \qquad l=1,\dots,m-1,
\\
\frac d{dz}y_m
=\frac{\sigma_1(\bbeta)-z\sigma_1(\balpha)}{z(z-1)}y_m
+\frac{\sigma_2(\bbeta)-z\sigma_2(\balpha)}{z(z-1)}y_{m-1}
+\dots+\frac{\sigma_m(\bbeta)-z\sigma_m(\balpha)}{z(z-1)}y_1.
\endgathered
\tag48
$$
Soglasno teoreme Liuvillya~\cite{17, gl.~3, \S\,27, p.~6}
on udovletvoryaet uravneniyu
$$
\frac d{dz}y=\Tr Q(z)\cdot y,
\tag49
$$
gde $\Tr Q(z)=(\beta-z\alpha)/(z(z-1))$~-- sled matricy
sistemy~\thetag{48}. Uravnenie~\thetag{49} mozhno perepisat\char"7E{}
v vide~\thetag{46}. Ego rexenie~\thetag{47} nakhodit\hbox{s}ya
neposredstvennym integrirovaniem.
\enddemo\cyrm

Sformuliruem nekotorye dostatoqnye usloviya iz~\cite{18}
na parametry funkcii \thetag{28} dlya togo, qtoby
funkcii~\thetag{29} byli algebraiqeski nezavi\-simy nad polem~$\CC(z)$:
\roster
\item"1)" {\cyri line\char"1Anaya neprivodimost\char"7E\/}:
$\alpha_l-\beta_j\notin\ZZ$ dlya vsekh $l,j=1,\dots,m$;
\item"2)" {\cyri neprivodimost\char"7E{} Belogo\/}~\cite{18, gl.~3, lemma~3.5.3}:
dlya lyubo\char"1A{} pary natural\char"7Enykh qisel~$m_1,m_2$, $m_1+m_2=m$,
ne suwestvuet qisel $u,v\in\QQ$ takikh, qto libo
$$
\gather
(\alpha_1,\alpha_2,\dots,\alpha_m)
\sim\biggl(\frac u{m_1},\frac{u+1}{m_1},\dots,\frac{u+m_1-1}{m_1},
\frac v{m_2},\frac{v+1}{m_2},\dots,\frac{v+m_2-1}{m_2}\biggr),
\\
(\beta_1,\beta_2,\dots,\beta_m)
\sim\biggl(\frac{u+v}m,\frac{u+v+1}m,\dots,\frac{u+v+m-1}m\biggr),
\endgather
$$
libo
$$
\gather
(\alpha_1,\alpha_2,\dots,\alpha_m)
\sim\biggl(\frac{u+v}m,\frac{u+v+1}m,\dots,\frac{u+v+m-1}m\biggr),
\\
(\beta_1,\beta_2,\dots,\beta_m)
\sim\biggl(\frac u{m_1},\frac{u+1}{m_1},\dots,\frac{u+m_1-1}{m_1},
\frac v{m_2},\frac{v+1}{m_2},\dots,\frac{v+m_2-1}{m_2}\biggr);
\endgather
$$
\item"3)" {\cyri kummerovskaya neprivodimost\char"7E\/}~\cite{18, gl.~3, lemma~3.5.6}:
ne suwestvuet delitelya~$m_0\ge2$ qisla~$m$ takogo, qto
$$
\gather
(\alpha_1,\alpha_2,\dots,\alpha_m)
\sim\biggl(\alpha_1+\frac1{m_0},\alpha_2+\frac1{m_0},
\dots,\alpha_m+\frac1{m_0}\biggr),
\\
(\beta_1,\beta_2,\dots,\beta_m)
\sim\biggl(\beta_1+\frac1{m_0},\beta_2+\frac1{m_0},
\dots,\beta_m+\frac1{m_0}\biggr);
\endgather
$$
\item"4)" $2\gamma\notin\ZZ$, gde
$\gamma=\alpha_1+\dots+\alpha_m-\beta_1-\dotsb-\beta_m$.
\endroster
Zapis\char"7E{}
$(\lambda_1,\dots,\lambda_m)\sim(\lambda_1',\dots,\lambda_m')$
oznaqaet, qto dlya nekotoro\char"1A{} perestanov\-ki
$\sigma\:\{1,\dots,m\}\to\{1,\dots,m\}$
pri vsekh $l=1,\dots,m$ vypolneno
$\lambda_l-\lambda_{\sigma(l)}'\in\ZZ$.

\proclaim{\cyrc\indent Lemma 16}\cyri
Pust\char"7E{} qisla
$\alpha_1,\dots,\alpha_m,\beta_1,\dots,\beta_m\in
\QQ\setminus\{-1,-2,\dots\}$
udovlet\-voryayut usloviyam~\rom{1)}--\rom{4)}.
Togda funkcii~\thetag{29}\hlop\rom, gde $f(z)$
opredelyaet\hbox{s}ya ryadom~\thetag{28}\hlop\rom,
algebraiqeski nezavisimy nad polem~$\CC(z)$.
\endproclaim\cyrm

\demo{\cyrc\indent Dokazatel\char"7Estvo}\cyrm
Soglasno~\cite{18, gl.~3, teorema~3.5.8}
v sluqae vypolneniya uslovi\char"1A~1)--4)
gruppa Galua odnorodnogo line\char"1Anogo
differencia\-l\char"7Enogo uravneniya~\thetag{45}
poryadka~$m$ izomorfna gruppe~$\operatorname{SL}_m(\CC)$.
E1to oznaqaet, qto funkcii, vkhodyawie v fundamental\char"7Enuyu
sistemu rexeni\char"1A{} uravneniya \thetag{45}, svyazany edinstvennym
algebraiqeskim sootnoxeniem nad polem $\CC(z)$ --
opredelitel\char"7E{} e1to\char"1A{}
fundamental\char"7Eno\char"1A{} sistemy rexeni\char"1A{}
yavlyaet\hbox{s}ya algeb\-raiqesko\char"1A{} funkcie\char"1A{} (lemma~15).
Sledovatel\char"7Eno, esli
$g(z)$~-- lyuboe netri\-vial\char"7Enoe rexenie uravneniya~\thetag{45},
to funkcii $g(z),\delta g(z),\dots,\delta^{m-1}g(z)$
algebraiqeski nezavisimy nad~$\CC(z)$.
Iz teoremy Nesterenko~\cite{19, teorema~2}
(sm\. takzhe~\cite{14, gl.~9, \S\,6, teorema~2})
sleduet, qto libo funkcii~\thetag{29} algebraiqeski
nezavisimy nad~$\CC(z)$, libo vse oni prinadlezhat
$\CC(z)$. Vtoro\char"1A{} sluqa\char"1A{} nevozmozhen, poskol\char"7Eku
funkciya~\thetag{28} ne yavlyaet\hbox{s}ya racional\char"7Eno\char"1A.
Lemma dokazana.
\enddemo\cyrm

\proclaim{\cyrc\indent Teorema 6}\cyri
Pust\char"7E{} parametry $\alpha_1,\dots,\alpha_m$ i $\beta_1,\dots,\beta_m$
funkcii~\thetag{28} udov\-letvoryayut usloviyam~\rom{1)}--\rom{4)}\rom,
$\beta_1,\dots,\beta_m$ poparno razliqny\rom,
$\xi$~-- nekotoroe raci\-onal\char"7Enoe qislo\rom,
$\xi=a_1/a_2\ne0$\rom, $a_2=\den\xi\in\NN$\rom, i
$\epsilon<1/(m+2)$~-- proizvol\char"7Enaya
polozhitel\char"7Enaya postoyannaya.
Oboznaqim qerez~$b_0$ naimen\char"7Exi\char"1A{}
obwi\char"1A{} znamenatel\char"7E{} qisel
$\gamma,\beta_1,\dots,\beta_m$\rom,
qerez $q_1$ i~$q_2$ proizvedeniya znamenatele\char"1A{} qisel
$\alpha_1,\dots,\alpha_m$ i $\beta_1,\dots,\beta_m$ sootvet\hbox{s}tvenno\rom,
qerez~$b$ naimen\char"7Exee obwee kratnoe qisel~$q_1,q_2$\rom,
qerez~$H$ maksimum modulya koe1fficientov mnogoqlenov~\thetag{31}.
Polozhim
$\Phi=e^{\rho(\den\beta_1)+\dots+\rho(\den\beta_m)}q_1/b$,
$$
\gathered
C_0=\bigl(8b_0He^{\chi(b_0)+3}\bigr)^{\epsilon(1-\log\epsilon)}
\Phi^{1+\epsilon+(2-(m-1)\epsilon)/(\epsilon^m(m-1)!)},
\\
\eta_0=\frac{(1+\epsilon)\log a_2+\log C_0}
{(1-(m+2)\epsilon)\log a_2-\log C_0-(2-(m+1)\epsilon)\log|a_1|},
\endgathered
$$
gde funkcii $\chi(\,\cdot\,)$ i $\rho(\,\cdot\,)$ opredelyayut\hbox{s}ya
formulami~\thetag{5} i~\thetag{44} sootvet\hbox{s}tvenno.
Esli dlya zadannogo $\xi$ vypolneno uslovie $\eta_0>0$\rom,
inymi slovami\rom, esli
$$
a_2^{1-(m+2)\epsilon}>C_0|a_1|^{2-(m+1)\epsilon},
$$
to qislo $f(\xi)$ irracional\char"7Eno. Bolee togo\rom,
dlya lyubogo $\eta>\eta_0$ i proizvol\char"7Enykh $p\in\ZZ$\rom, $q\in\NN$\rom,
$q>q_*(\xi,\epsilon,\eta)$\rom,
spravedliva ocenka
$$
\biggl|f(\xi)-\frac pq\biggr|>q^{-1-\eta}.
$$
\endproclaim\cyrm

\demo{\cyrc\indent Dokazatel\char"7Estvo}\cyrm
Pust\char"7E{} $T_0(z)=bT(z)=bz(z-1)$~--
obwi\char"1A{} znamenatel\char"7E{} sistemy differencial\char"7Enykh
uravneni\char"1A~\thetag{32},
kotoro\char"1A{} udovletvoryayut sovokupnost\char"7E{}
funkci\char"1A~\thetag{29} iz klassa $\GG(1,\Phi)$
(sm\. sledstvie iz lemmy~14);
pri e1tom
$$
T_0(z)Q_{lj}(z)\in\ZZ[z],
\qquad l=1,\dots,m, \quad j=0,\dots,m.
$$
Togda
$$
\max\Bigl\{\deg T_0-1,\max_{l,j}\{\deg T_0Q_{lj}\}\Bigr\}=1,
\qquad
\max\Bigl\{H(T_0),\max_{l,j}\{H(T_0Q_{lj})\}\Bigr\}=bH.
$$
Iz usloviya~4) sleduet, qto $\gamma\ne0$.
Soglasno lemme~11 sokrawenie faktorialov dlya sistemy differencial\char"7Enykh
uravneni\char"1A~\thetag{33}, sopryazhenno\char"1A{} k odnorodno\char"1A{} qasti sistemy~\thetag{32},
proiskhodit s postoyanno\char"1A{} $\Psi=b_0e^{\chi(b_0)+2}$
(s postoyanno\char"1A{} $\Psi=b_0e^{\chi(b_0)+2}/b$, esli v opredelenii~3
mnogoqlen~$T(z)$ zamenit\char"7E{} na~$T_0(z)$).
Iz lemmy~16 sleduet algebraiqeskaya nezavisimost\char"7E{} funkci\char"1A~\thetag{29}
nad polem~$\CC(z)$. Primenyaya teper\char"7E{} osnovnuyu teoremu
i nera\-venstva~(0.9) iz raboty~\cite{16}, poluqaem trebuemoe
utverzhdenie.
\enddemo\cyrm

\subhead\cyrb
7. Sokrawenie faktorialov dlya line\char"1Anykh odnorodnykh
differencial\char"7Enykh uravneni\char"1A{} s postoyannymi koe1fficientami
\endsubhead\cyrm
Rassmotrim
\linebreak
sistemu line\char"1Anykh odnorodnykh
differencial\char"7Enykh uravneni\char"1A{}
$$
\frac d{dz}y_l=\sum_{j=1}^mA_{lj}y_j,
\quad l=1,\dots,m,
\qquad A_{lj}\in\CC, \quad l,j=1,\dots,m,
\tag50
$$
i svyazanny\char"1A{} s ne\char"1A{} differencial\char"7Eny\char"1A{} operator
$$
\sA=[A]=\sum_{l=1}^m\sum_{j=1}^mA_{lj}y_j\frac\partial{\partial y_l},
\qquad
A=\pmatrix
A_{11} & \hdots & A_{1m} \\
\hdotsfor3 \\
A_{m1} & \hdots & A_{mm}
\endpmatrix.
\tag51
$$
Takim obrazom, s kazhdo\char"1A{} matrice\char"1A~$A$
mozhno svyazat\char"7E{} differencial\char"7Eny\char"1A{} operator~$[A]$
soglasno formule~\thetag{51}.
Otmetim, prezhde vsego, svo\char"1Astvo line\char"1Anosti
$$
[\lambda_1A_1+\lambda_2A_2]
=\lambda_1[A_1]+\lambda_2[A_2],
\qquad \lambda_1,\lambda_2\in\CC,
$$
dlya lyubykh kvadratnykh matric $A_1$ i~$A_2$ razmera~$m$.

Dlya differencial\char"7Enykh operatorov~$[A]$ i~$[B]$
vvedem operacii formal\char"7E\-nogo umnozheniya
$$
\align
[A]\cdot[B]
&=\biggl(\sum_{l=1}^m\sum_{j=1}^mA_{lj}y_j\frac\partial{\partial y_l}\biggr)
\cdot\biggl(\sum_{i=1}^m\sum_{k=1}^mB_{ik}y_k\frac\partial{\partial y_i}\biggr)
\\
&=\sum_{l=1}^m\sum_{i=1}^m
\biggl(\sum_{j=1}^mA_{lj}y_j\biggr)
\biggl(\sum_{k=1}^mB_{ik}y_k\biggr)
\frac{\partial^2}{\partial y_l\,\partial y_i}
=[B]\cdot[A]
\endalign
$$
i kompozicii
$$
\align
[B]\circ[A]
&=\biggl(\sum_{i=1}^m\sum_{k=1}^mB_{ik}y_k\frac\partial{\partial y_i}\biggr)
\circ\biggl(\sum_{l=1}^m\sum_{j=1}^mA_{lj}y_j\frac\partial{\partial y_l}\biggr)
\\
&=\sum_{i=1}^m\sum_{k=1}^mB_{ik}y_k\sum_{l=1}^mA_{li}
\frac\partial{\partial y_l}
+\biggl(\sum_{i=1}^m\sum_{k=1}^mB_{ik}y_k\frac\partial{\partial y_i}\biggr)
\cdot\biggl(\sum_{l=1}^m\sum_{j=1}^mA_{lj}y_j\frac\partial{\partial y_l}\biggr)
\\
&=\sum_{l=1}^m\sum_{k=1}^m\sum_{i=1}^m
A_{li}B_{ik}y_k\frac\partial{\partial y_l}
+\sum_{l=1}^m\sum_{i=1}^m
\biggl(\sum_{j=1}^mA_{lj}y_j\biggr)
\biggl(\sum_{k=1}^mB_{ik}y_k\biggr)
\frac{\partial^2}{\partial y_l\,\partial y_i}
\\
&=[AB]+[A]\cdot[B].
\endalign
$$
Pomimo oqevidno\char"1A{} associativnosti kazhdaya iz e1tikh
operaci\char"1A{} distributivna otnositel\char"7Eno slozheniya.
Kommutativnost\char"7E{} imeet mesto tol\char"7Eko dlya operacii
formal\char"7Enogo umnozheniya, khotya poryadok mnozhitele\char"1A{}
v kompozi\-cii
$$
\sA_n=([A]-n+1)\circ([A]-n+2)\circ\dotsb\circ([A]-1)\circ[A]
\tag52
$$
ne imeet znaqeniya.
Otmetim takzhe ``differencial\char"7Enoe
svo\char"1Astvo'' ope\-racii kompozicii:
$$
[B]\circ([A_1]\cdot[A_2])
=([B]\circ[A_1])\cdot[A_2]
+[A_1]\cdot([B]\circ[A_2]).
$$
Pod $[A]^n$ budem ponimat\char"7E{} formal\char"7Enoe proizvedenie
$n$~operatorov~$[A]$.

\definition{\cyrc\indent Opredelenie 4}\cyrm
Skazhem, qto {\cyri differencial\char"7Eny\char"1A{} operator\/}~\thetag{51}
(ili {\cyri sistema uravneni\char"1A\/}~\thetag{50}),
otveqayuwi\char"1A{} matrice~$A$, {\cyri udovletvoryaet
usloviyu sokraweniya faktorialov s postoyanno\char"1A\/} $\Psi\ge1$,
esli suwestvuet posledovatel\char"7Enost\char"7E{} natural\char"7Enykh qisel
$\{\psi_k\}_{k\in\NN}$ takaya, qto operatory
$$
\psi_k\frac1{n!}\sA_n,
\qquad n=0,1,\dots,k, \quad k\in\NN,
\tag53
$$
perevodyat kol\char"7Eco $\ZZ[y_1,\dots,y_m]$ v sebya i
$$
\limsup_{k\to\infty}\psi_k^{1/k}\le\Psi.
$$
\enddefinition\cyrm

\proclaim{\cyrc\indent Lemma 17}\cyri
Dlya differencial\char"7Enogo operatora~\thetag{52}
vypolneno tozhdest\-vo
$$
\frac1{n!}\sA_n
=\sum\Sb s_1,s_2,\dots,s_n\ge0\\s_1+2s_2+\dots+ns_n=n\endSb
\frac1{s_1!}[A]^{s_1}\cdot
\frac1{s_2!}\biggl[\frac{\<A\>_2}{2!}\biggr]^{s_2}
\dotsb\frac1{s_n!}\biggl[\frac{\<A\>_n}{n!}\biggr]^{s_n},
\tag54
$$
gde simvol $\<\,\cdot\,\>_n$ opredelyaet\hbox{s}ya formulo\char"1A~\thetag{1}.
\endproclaim\cyrm

\demo\nofrills{\cyrc\indent Dokazatel\char"7Estvo}\cyrm\
provedem indukcie\char"1A{} po~$n$. Baza indukcii $n=1$
oqevidna. Pust\char"7E{} tozhdestvo~\thetag{54} spravedlivo
dlya zadannogo~$n$. Polozhim
$$
U(s_1,s_2,\dots,s_n)
=\frac1{s_1!}[A]^{s_1}\cdot
\frac1{s_2!}\biggl[\frac{\<A\>_2}{2!}\biggr]^{s_2}
\dotsb\frac1{s_n!}\biggl[\frac{\<A\>_n}{n!}\biggr]^{s_n}
\tag55
$$
i perepixem~\thetag{54} v vide
$$
\align
\frac1{n!}\sA_n
&=\sum\Sb s_1,s_2,\dots,s_n\ge0\\s_1+2s_2+\dots+ns_n=n\endSb
U(s_1,s_2,\dots,s_n)
\\
&=\sum\Sb s_1,s_2,\dots,s_n,s_{n+1}\ge0\\
s_1+2s_2+\dots+ns_n+(n+1)s_{n+1}=n\endSb
U(s_1,s_2,\dots,s_n,s_{n+1})
\endalign
$$
(zdes\char"7E, oqevidno, $s_{n+1}=0$).
Poskol\char"7Eku
$$
([A]-i)\circ[\<A\>_i]=[A]\cdot[\<A\>_i]+[\<A\>_{i+1}],
$$
imeem
$$
\gather
\aligned
&
([A]-n)\circ U(s_1,s_2,\dots,s_n,s_{n+1})
\\ &\qquad
=(s_1+1)U(s_1+1,s_2,\dots,s_n,s_{n+1})
\\ &\qquad\quad
+2(s_2+1)U(s_1-1,s_2+1,s_3,\dots,s_n,s_{n+1})
\\ &\qquad\quad
+3(s_3+1)U(s_1,s_2-1,s_3+1,s_4,\dots,s_n,s_{n+1})
+\dotsb
\\ &\qquad\quad
+(n+1)(s_{n+1}+1)U(s_1,s_2,\dots,s_{n-1},s_n-1,s_{n+1}+1),
\endaligned
\\
s_1,s_2,\dots,s_n,s_{n+1}\ge0, \qquad
s_1+2s_2+\dots+ns_n+(n+1)s_{n+1}=n.
\endgather
$$
Poe1tomu
$$
\align
&
\frac1{(n+1)!}\sA_{n+1}
=\frac{[A]-n}{n+1}\circ\frac1{n!}\sA_n
\\ &\qquad
=\sum\Sb s_1,s_2,\dots,s_n,s_{n+1}\ge0\\
s_1+2s_2+\dots+ns_n+(n+1)s_{n+1}=n\endSb
\frac{s_1+1}{n+1}U(s_1+1,s_2,\dots,s_n,s_{n+1})
\\ &\qquad\quad
+\frac{2(s_2+1)}{n+1}U(s_1-1,s_2+1,s_3,\dots,s_n,s_{n+1})
+\dotsb
\\ &\qquad\quad
+\frac{(n+1)(s_{n+1}+1)}{n+1}U(s_1,s_2,\dots,s_{n-1},s_n-1,s_{n+1}+1)
\\ &\qquad
=\sum\Sb s_1',s_2',\dots,s_{n+1}'\ge0\\
s_1'+2s_2'+\dots+(n+1)s_{n+1}'=n+1\endSb
\frac{s_1'+2s_2'+\dots+(n+1)s_{n+1}'}{n+1}
U(s_1',s_2',\dots,s_{n+1}')
\\ &\qquad
=\sum\Sb s_1',s_2',\dots,s_{n+1}'\ge0\\
s_1'+2s_2'+\dots+(n+1)s_{n+1}'=n+1\endSb
U(s_1',s_2',\dots,s_{n+1}'),
\endalign
$$
qto oznaqaet spravedlivost\char"7E{} tozhdestva~\thetag{54} dlya~$n+1$.
Lemma dokazana.
\enddemo\cyrm

\proclaim{\cyrc\indent Lemma 18}\cyri
Esli koe1fficienty matricy~$B$ celoqislenny\rom, to
differencial\char"7Eny\char"1A{} operator $[B]^s/s!$
perevodit kol\char"7Eco $\ZZ[y_1,\dots,y_m]$ v sebya.
\endproclaim\cyrm

\demo{\cyrc\indent Dokazatel\char"7Estvo}\cyrm
Zapixem differencial\char"7Eny\char"1A{} operator~$[B]$ v vide
$$
[B]=\sum_{l=1}^mb_l\frac\partial{\partial y_l},
\qquad b_l=\sum_{j=1}^mB_{lj}y_j\in\ZZ[y_1,\dots,y_m].
$$
Togda differencial\char"7Eny\char"1A{} operator
$$
\frac1{s!}[B]^s
=\frac1{s!}\biggl(\sum_{l=1}^mb_l\frac\partial{\partial y_l}\biggr)^s
=\sum\Sb l_1,\dots,l_m\ge0\\l_1+\dots+l_m=s\endSb
\frac{b_1^{l_1}}{l_1!}\frac{\partial^{l_1}}{\partial y_1^{l_1}}
\dotsb
\frac{b_m^{l_m}}{l_m!}\frac{\partial^{l_m}}{\partial y_m^{l_m}}
$$
perevodit kol\char"7Eco $\ZZ[y_1,\dots,y_m]$ v sebya,
tak kak e1tim svo\char"1Astvom obladayut operatory
$$
\frac1{l_j!}\frac{\partial^{l_j}}{\partial y_j^{l_j}},
\qquad j=1,\dots,m
$$
(primer~1), i, krome togo,
$$
b_j^{l_j}\in\ZZ[y_1,\dots,y_m],
\qquad j=1,\dots,m.
$$
Lemma dokazana.
\enddemo\cyrm

\proclaim{\cyrc\indent Lemma 19}\cyri
Pust\char"7E{} $s_1,s_2,\dots,s_k$~-- celye neotricatel\char"7Enye qisla\rom,
$p$~-- prostoe qislo. Togda
$$
s_1\tau_p(1)+s_2\tau_p(2)+\dots+s_k\tau_p(k)
\le\tau_p(s_1+2s_2+\dots+ks_k),
\tag56
$$
gde $\tau_p(k)$~-- stepen\char"7E{} vkhozhdeniya qisla~$p$ v~$k!$
\rom(sm.~\thetag{4}\rom).
\endproclaim\cyrm

\demo{\cyrc\indent Dokazatel\char"7Estvo}\cyrm
Otmetim, prezhde vsego, qto dlya de\char"1Astvitel\char"7Enykh
$\xi_1,\allowmathbreak\dots,\xi_k$ imeet mesto neravenstvo
$$
[\xi_1]+\dots+[\xi_k]\le[\xi_1+\dots+\xi_k],
$$
gde $[\,\cdot\,]$~-- celaya qast\char"7E{} qisla
(sm.~\cite{2, otd.~8, gl.~1, \S\,1, zadaqa~7}). Ot\hbox{s}yuda sleduet,
qto
$$
s_1[\xi_1]+s_2[\xi_2]+\dots+s_k[\xi_k]
\le[s_1\xi_1+s_2\xi_2+\dots+s_k\xi_k].
\tag57
$$
Posledovatel\char"7Eno podstavlyaya v~\thetag{57} znaqeniya
$$
\xi_n=\frac np,\frac n{p^2},\frac n{p^3},\dots,
\qquad n=1,\dots,k,
$$
i summiruya poluqennye neravenstva,
soglasno formule~\thetag{4} poluqaem~\thetag{56}.
Lemma dokazana.
\enddemo\cyrm

\proclaim{\cyrc\indent Teorema 7}\cyri
Pust\char"7E{} differencial\char"7Eny\char"1A{} operator~\thetag{51}
otveqaet racional\char"7E\-no\char"1A{} matrice~$A$\rom, $\den A=b$\rom,
minimal\char"7Eny\char"1A{} mnogoqlen kotoro\char"1A{}
ne imeet kratnykh korne\char"1A\rom;
$t_1,t_2$~-- naimen\char"7Exie obwie znamenateli e1lementov
matric~$T$ i $T^{-1}$ sootvet\hbox{s}tvenno\rom,
gde $T$~-- matrica perekhoda ot~$A$ k ee
zhordanovo\char"1A{} normal\char"7Eno\char"1A{} forme.
Togda operator~$[A]$ udovletvoryaet usloviyu sokraweniya
faktorialov s postoyanno\char"1A~$t_1t_2be^{\chi(b)}$.
\endproclaim\cyrm

\demo{\cyrc\indent Dokazatel\char"7Estvo}\cyrm
Pust\char"7E{} $k$~-- proizvol\char"7Enoe natural\char"7Enoe qislo.
Soglasno lemme~9 matricy
$$
t_1t_2b^n\prod_{p\mid b}p^{\tau_p(n)}\cdot\frac{\<A\>_n}{n!},
\qquad n=0,1,\dots,k,
$$
imeyut celoqislennye e1lementy.
Poe1tomu iz lemmy~18 sleduet, qto
differencial\char"7Eny\char"1A{} operator~\thetag{55} posle umnozheniya na
$$
(t_1t_2)^nb^n\prod_{p\mid b}p^{s_1\tau_p(1)+s_2\tau_p(2)+\dots+s_n\tau_p(n)}
\tag58
$$
perevodit kol\char"7Eco $\ZZ[y_1,\dots,y_m]$ v sebya. Poskol\char"7Eku
$s_1+2s_2+\dots+ns_n=n\le k$, iz lemmy~19 sleduet, qto qislo~\thetag{58}
delit
$$
\psi_k=(t_1t_2)^kb^k\prod_{p\mid b}p^{\tau_p(k)}.
$$
Primenim tozhdestvo lemmy~17. Dlya posledovatel\char"7Enosti
$\{\psi_k\}_{k\in\NN}$ opera\-tory~\thetag{53}
perevodyat kol\char"7Eco $\ZZ[y_1,\dots,y_m]$ v sebya.
Ocenka~\thetag{18} zaverxaet dokazatel\char"7Estvo teoremy.
\enddemo\cyrm

\remark{\cyrc\indent Zameqanie}\cyrm
S sistemo\char"1A{} differencial\char"7Enykh
uravneni\char"1A~\thetag{8},~\thetag{12}
fuksov\-skogo tipa svyazan differencial\char"7Eny\char"1A{} operator
$$
D=\frac d{dz}+\frac1{z-\gamma_1}[A_1]+\dots
+\frac1{z-\gamma_s}[A_s].
$$
E1to obstoyatel\char"7Estvo pozvolyaet s pomow\char"7Eyu teoremy~7 dat\char"7E{}
drugoe dokazatel\char"7Estvo teoremy~5 v sluqae poparno kommutiruyuwikh
matric $A_1,\dots,A_s$, minimal\char"7Eny\char"1A{} mnogoqlen kazhdo\char"1A{} iz kotorykh
ne imeet kratnykh korne\char"1A.
Pri e1tom, odnako, poluqaet\hbox{s}ya neskol\char"7Eko khudxaya (po sravneniyu
s teoremo\char"1A~5) postoyannaya~$\Psi$.
\endremark\cyrm

\addto\eightpoint{\normalbaselineskip=.95\normalbaselineskip\normalbaselines}

\enddocument
\end